\documentclass[arXiv,noinfoline]{imsart}

\usepackage{xr}
\usepackage{amsmath,amstext,amssymb}
\arxiv{arXiv:0000.0000}

\usepackage{xparse}
\usepackage{lipsum} 
\usepackage{algorithm}
\usepackage{algpseudocode}
\usepackage{accents} 
\usepackage{etoolbox}

\usepackage{amsbsy} 
\usepackage{amsthm}
\usepackage{mathrsfs}
\usepackage{eufrak}
\usepackage{mathrsfs}
\usepackage{color}
\usepackage{verbatim}
\usepackage{graphicx}
\usepackage{bm}
\usepackage{enumerate}
\usepackage{epstopdf} 
\usepackage{natbib}
\RequirePackage[colorlinks,citecolor=blue,urlcolor=blue]{hyperref}
\usepackage{subfig}
\usepackage[final]{pdfpages}

\usepackage{algorithm}  
\usepackage{algpseudocode} 
\usepackage{algorithmicx}     

\usepackage{soul}  

\usepackage[dvipsnames]{xcolor}

\newcommand{\bx}{\bm{x}}

\newcommand{\HaarMu}{\mu}

\newcommand{\ConstOne}{K}

\usepackage{xstring}
\usepackage[normalem]{ulem}
\definecolor{ultramarine}{RGB}{38,29,163}

\newcommand{\qedwhite}{\hfill \ensuremath{\Box}}
\newcommand{\SpaceD}{\mathcal{S}_p}

\newcommand{\SpaceM}{\widetilde{\mathcal{V}}_{n,p}}
\newcommand{\SpaceV}{\mathcal{V}_{p,p}}

\newcommand{\StiefelS}{\mathcal{V}_{n,p}}

\newcommand{\ML}{{\cal{ML}}}

\newcommand{\bbeta}{\bm{\beta}}

\newcommand{\bd}{\bm{d}}
\newcommand{\BoEta}{\boldsymbol{\eta}}
\newtheorem{rmrk}{Remark}
\newtheorem{theorem}{Theorem}
\newtheorem{lemma}{Lemma}
\newtheorem{defn}{Definition}

\newcommand{\priorXzero}{\Psi}
\newcommand{\iMat}{\mathbf{I}_{p}}

\newcommand{\normtwo}[1]{{\left\lVert#1\right\rVert}_2}

\newcommand{\hyp}{{}_0F_1\left(\frac{n}{2},\frac{D^2}{4}\right)}
\newcommand{\hypinline}{{}_0F_1\left({n}/{2},{D^2}/{4}\right)}
\newcommand{\hypinlineF}{{}_0F_1\left({n}/{2},{F^TF}/{4}\right)}
\newcommand{\partialhyp}[1]{\frac{\partial}{\partial\,{d_{#1}}}\,\left[\hyp\right]}

\newcommand{\IMDY}{{\it{CCPC}}}
\newcommand{\JMDY}{{\it{JCPC}}}

\newcommand{\CCPD}{{\it{CCPD}}}
\newcommand{\JCPD}{{\it{JCPD}}}

\definecolor{auburn}{rgb}{0.53, 0.1, 0.5}

\newcommand{\MLDensity}{f_{\ML}}
\newcommand{\posterior}

\newcommand{\variableX}{\bd}

\newcommand{\IndVzero}[1]{\mathbb{I}({X\in \mathcal{V}^{#1}_0})}
\newcommand{\Rnp}{\mathbb{R}^{n \times p}}
\newcommand{\Rpp}{\mathbb{R}^{p \times p}}
\newcommand{\vecnorm}[1]{\lVert #1\rVert}

\newcommand{\Rplus}{\mathbb{R}_{+}}

\newcommand{\m}{{\bf m_{\BoEta}}} 

\newcommand{\SetWithMode}{\mathcal{S}}
\newcommand{\SetWithModePrime}{\mathcal{S}}

\newcommand{\hyparam}[2]{
    \IfEqCase{#1}{
        {M}{\xi^{#2}_c}
        {V}{\gamma^{#2}_c}%
        
    }
  }
  
  \newcommand{\hyparamM}[1]{\xi^{#1}}
   \newcommand{\hyparamV}[1]{\gamma^{#1}}

\ExplSyntaxOn
\DeclareExpandableDocumentCommand{\Measure}{o m}
 {
 \ifstrempty{#1}{
  \str_case:nn { #2 }
   {
    {M}{\mu}%
    {D}{\mu\sb{1}}%
    {V}{\mu\sb{2}}
    {X}{\mu}
    {B}{\mu}
   }
 }{
   \str_case:nn { #1 }
   {
    {1}{
    \str_case:nn { #2 }
   {
     {M}{d\mu(M)}%
        {D}{d\mu\sb{1}(\bd)}%
        {V}{d\mu\sb{2}(V)}
        {X}{d\mu(X)}
        {Y}{d\mu(Y)}
        {MDV} {d\mu(M)\; d\mu\sb{1}(\bd) \;d\mu\sb{2}(V) }
   }
    }%
     {2}{
    \str_case:nn { #2 }
   {
    {M}{d\mu(M^{\ast})}%
        {D}{d\mu\sb{1}(\bd^{\ast})}%
        {V}{d\mu\sb{2}(V^{\ast})}
        {X}{d\mu(X^{\ast})}
   }
    }
    {3}{
    \str_case:nn { #2 }
   {
     {M}{\mu(dM^{\star})}%
        {D}{\mu_2(d\bd^{\star})}%
        {V}{\mu_1(dV^{\star})}
        {X}{\mu(X^{\star})}
   }
    }
   }

    }
 }
\ExplSyntaxOff

\newcommand{\fourpartdef}[8]
{
	\left\{
		\begin{array}{lll}
			#1 & \mbox{if } #2 \\
			#3 & \mbox{if } #4 \\
			#5 & \mbox{if } #6\\
			#7 & \mbox{if } #8
		\end{array}
	\right.
}

\newcommand{\AllData}{{\{W_i\}_{i=1}^{N}}}
\newtheorem{example}{Example}

\newlength{\dhatheight}
\newcommand{\doublehat}[1]{%
    \settoheight{\dhatheight}{\ensuremath{\widehat{#1}}}%
    \addtolength{\dhatheight}{-0.35ex}%
    \widehat{\vphantom{\rule{2pt}{\dhatheight}}%
    \smash{\hspace{-0.5mm}\widehat{#1}}}}
    
\newcommand{\LevelSetMDV}{\mathcal{A}}
\newcommand{\SetWithModeMDV}{\mathcal{A}}

\newcommand{\rr}{d}
\newcommand{\R}{D}

\newcommand{\TermK}{A}
\newcommand{\M}{M}
\newcommand{\p}{c}    
    
\newcommand{\yy}{\color{black}}
\newcommand{\jj}{\color{black} \rm}
\begin{document}

\begin{frontmatter}

\title{{ A Bayesian Approach for Analyzing Data on the Stiefel Manifold }}
\begin{aug}
\author{\fnms{Subhadip} \snm{Pal}\thanksref{m2}
\ead[label=e2]{subhadip.pal@louisville.edu}},
\author{\fnms{Subhajit} \snm{Sengupta}\thanksref{m1}\ead[label=e1]{subhajit@uchicago.edu}},
\author{\fnms{Riten} \snm{Mitra}\thanksref{m2}
\ead[label=e3]{r0mitr01@louisville.edu}} and
\author{\fnms{Arunava} \snm{Banerjee}\thanksref{t2,m4}\ead[label=e5]{ arunava@ufl.edu}}

\runauthor{Pal et. al.}

\address{Department of Bioinformatics and Biostatistics, University of Louisville\thanksmark{m2}}
\address{Center for Psychiatric Genetics, NorthShore University HealthSystem\thanksmark{m1}}
\address{Department of Computer \& Information Science \& Engineering, University of Florida\thanksmark{m4}}


\end{aug}

\begin{abstract}

Directional data emerges in a wide array of applications, ranging from atmospheric sciences to medical imaging. Modeling such data, however, poses unique challenges by virtue of their being constrained to non-Euclidean spaces like manifolds. Here, we present a unified Bayesian framework for inference on the Stiefel manifold using the Matrix Langevin distribution. Specifically, we propose a novel family of conjugate priors and establish a number of theoretical properties relevant to statistical inference. 
Conjugacy enables translation of these properties to their corresponding posteriors, which we exploit to develop the posterior inference scheme. For the implementation of the posterior computation, including the posterior sampling, we adopt a novel computational procedure for evaluating the hypergeometric function of matrix arguments that appears as normalization constants in the relevant densities.
\end{abstract}

\medskip



\end{frontmatter}
\newpage
\section{Introduction}
\label{sec:intro}

Analysis of directional data is a major area of investigation in statistics. Directional data range from unit vectors in the simplest case to sets of ordered orthonormal frames in the general scenario. Since the associated sample space is non-Euclidean, standard statistical methods developed for the Euclidean space may not be appropriate to analyze such data. Additionally, it is often desirable to design statistical methods that take into consideration the underlying geometric structure of the sample space. There is a need for methodological development for a general sample space such as the Stiefel manifold~\citep{James:1976,Chikuse:2012} that goes beyond those techniques designed for simpler non-Euclidean spaces like the circle or the sphere. Such a novel methodology can support various emerging applications, increasingly seen in the fields of Biology~\citep{Downs:1972,Mardia:1977}, Computer science~\citep{Turaga:2008,Lui:2008} and Astronomy~\citep{Mardia:2009,Lin:2017}, to mention but a few.




One of the most widely used probability distributions on the Stiefel manifold is the matrix Langevin distribution introduced by~\cite{Downs:1972}, also known as the Von-Mises Fisher matrix distribution~\citep{Mardia:2009,Khatri:1977}. In early work,~\cite{Mardia:1977} and~\cite{Jupp:1980} investigated properties of the matrix Langevin distribution and developed inference procedures in the  frequentist setup~\citep{Chikuse:2012}. The form of the maximum likelihood estimators and the profile likelihood estimators for the related parameters can be found in ~\cite{Jupp:1979, Mardia:1977, Chikuse:1991,Chikuse:1991:as, Chikuse:1998}. It is not patently clear from these works whether the form of the associated asymptotic variance can be obtained directly without using bootstrap procedures. A major obstacle facing the development of efficient inference techniques for this family of distributions has been the intractability of the corresponding normalizing constant, a hypergeometric function of a matrix argument~\citep{ Mardia:2009, Muirhead:2009,Gross:1989}. Inference procedures have been developed exploiting approximations that are available when the argument to this function is either small or large.
Almost all the hypothesis testing procedures~\citep{Jupp:1979, Mardia:1977, Chikuse:1991,Chikuse:1991:as, Chikuse:1998} therefore depend not only on large sample asymptotic distributions but also on the specific cases when the  concentration parameter is either large or small~\citep{Chikuse:2012,Mardia:1977,Downs:1972}. In particular, a general one sample or two sample hypothesis testing method for the finite sample case is yet to be developed.  

For any given dataset, the stipulation of large sample is comparatively easier to verify than checking whether the magnitude of the concentration is large. It may not be possible to ascertain whether the concentration is large  before the parameter estimation procedure, which is then confounded by the fact that the existing parameter estimation procedures themselves require the assumption of large concentration to work correctly. Hence, from a practitioner's point of view, it is often difficult to identify whether the above-mentioned procedures are suitable for use on a particular dataset.

Although a couple of Bayesian procedures have been proposed in related fields (see references in~\cite{Lin:2017}), a comprehensive Bayesian analysis is yet to be developed for the matrix Langevin distribution. In a recent paper,~\cite{Lin:2017} have developed a Bayesian mixture model of matrix Langevin distributions for clustering on the Stiefel manifold, where they have used a prior structure that does not have conjugacy. To accomplish posterior inference,~\cite{Lin:2017} have used a nontrivial data augmentation strategy based on a rejection sampling technique laid out in~\cite{Rao:2016}. It is worthwhile to note that the specific type of data augmentation has been introduced to tackle the intractability of the hypergeometric function of a matrix argument. 
It is well known that data augmentation procedures often suffer from slow rate of convergence~\citep{VanDyk:meng:2001,Hobert:2011}, particularly when combined with an inefficient rejection sampler. Elsewhere,~\cite{Hornik:2014} have proposed a class of conjugate priors but have not presented an inference procedure for the resulting posterior distributions.

In this article, we develop a comprehensive Bayesian framework for the matrix Langevin distribution, starting with the construction of a flexible class of conjugate priors, and proceeding all the way to the design of an efficient posterior computation procedure. The difficulties arising from the intractability of the normalizing constant do not, of course, disappear with the mere adoption of a Bayesian approach. We employ nontrivial strategies to derive a unique posterior inference scheme in order to handle the intractability of the normalizing constant. A key step in the proposed posterior computation is the evaluation of the hyper-geometric function of a matrix argument, that can be computed using the algorithm developed in \cite{Koev:2006}. Although general, this algorithm has certain limitations vis-\`a-vis measuring the precision of its output. We therefore construct a reliable and computationally efficient procedure to compute a specific case of the hypergeometric function of matrix argument, that has theoretical precision guarantees (Section~\ref{sec:HYPComputation}).
The procedure is applicable to a broad class of datasets including most, if not all, of the applications found in~\cite{Downs:1971, Downs:1972, Jupp:1979,Jupp:1980, Mardia:1977, Mardia:2007, Mardia:2009, Chikuse:1991:as,Chikuse:1991,Chikuse:1998,Chikuse:2003,Sei:2013, Lin:2017}. The theoretical framework proposed in this article is applicable to all matrix arguments regardless of dimensionality. 
In the following two paragraphs, we summarize our contributions.

We begin by adopting a suitable representation of the hypergeometric function of a matrix argument to view it as a function of a vector argument. We explore several of its properties that are useful for subsequent theoretical development, and also adopt an alternative parametrization of the matrix Langevin distribution so that the modified representation of the hypergeometric function can be used. When viewed as an exponential family of distributions, the new parameters of the matrix Langevin distribution are not the natural parameters ~\citep{Casella:2002}. Thus the construction of the conjugate prior does not directly follow from~\cite{Diaconis:Ylvisaker:1979} (DY), an issue that we elaborate on (Section~\ref{subsec:DY_inapplicable}). We then propose two novel and  reasonably large classes of conjugate priors, and based on theoretical properties of the matrix Langevin distribution and the hypergeometric function, we  establish their propriety. We study useful properties of the constructed class of distributions to demonstrate that the hyperparameters related to the class of distributions have natural interpretations. Specifically, the class of constructed distributions is characterized by two hyperparameters, one controls the \emph{location} of the distribution while the other determines the \emph{scale}. This interpretation not only helps us understand the nature of the class of distributions but also aids in the selection of hyperparameter settings. The constructed class of  prior distributions is flexible because one can incorporate prior knowledge via appropriate hyperparameter selection; and at the same time, in the absence of prior knowledge, there is a provision to specify the hyperparameters to construct a uniform prior. Since this uniform prior is improper by nature, we extend our investigation to identify the conditions under which the resulting posterior is a proper probability distribution.


Following this, we discuss properties of the posterior and inference. We show unimodality of the resulting posterior distributions and derive a computationally efficient expression for the posterior mode. We also demonstrate that the posterior mode is a consistent estimator of the related parameters. We develop a Gibbs sampling algorithm  to sample from the resulting posterior distribution. One of the conditionals in the Gibbs sampling algorithm is a novel class of distributions that we have introduced in this article for the first time. We develop and make use of properties such as unimodality and log-concavity to derive an efficient rejection sampler to sample from this distribution. We perform multiple simulations to showcase the generic nature of our framework and to report estimation efficiency for the different algorithms. We end with an application demonstrating the strength of our approach.

 
We should note that a significant portion of the article is devoted to establishing a number of novel properties of the hypergeometric function of matrix arguments. These properties play a key role in the rigorous development of the statistical procedures. These properties, including the exponential type upper and lower bounds for the function, may also be relevant to a broader range of scientific disciplines.



The remainder of the article is organized as follows. In Section~\ref{sec:stiefel_distr}, we introduce the matrix Langevin distribution defined on the Stiefel manifold and explore some of its important properties. Section~\ref{sec:prior_construct} begins with a discussion of the inapplicability of DY's theorem, following which we present the construction of the conjugate prior for the parameters of the matrix Langevin distribution. In particular, we establish propriety of a class of posterior and prior distributions by proving the finiteness of the integral of  specific  density kernels. In Section~\ref{sec:hyperparameter_singleML} and~\ref{sec:posteriorProperties}, we lay out the hyperparameter selection procedure and derive properties of the posterior. In Section~\ref{subsubsec:post_comp} we develop the posterior inference scheme. In Sections~\ref{sec:simul_data_app} and~\ref{sec:real_data_app}, we validate the robustness of our framework with experiments using simulated datasets and demonstrate the applicability of the framework using a real dataset, respectively. Finally, in Section~\ref{sec:disc}, we discuss other developments and a few possible directions for future research. Proofs of all theorems and properties of the hypergeometric function of matrix arguments are deferred to the supplementary material.



\paragraph{Notational Convention}
\begin{itemize}
\renewcommand\labelitemi{}
\item $\mathbb{R}^p$ = The $p$-dimensional Euclidean space.
\item $\mathbb{R}_{+}^{p} = \left\{\left(x_1, \ldots, x_p\right) \in \mathbb{R}^{p} : 0< x_i \text{ for } i =1,\ldots p  \right\}$.
\item $\mathcal{S}_p= \left\{\left(d_1, \ldots, d_p\right) \in \mathbb{R}_{+}^{p} : 0< d_p< \cdots < d_1 <\infty \right\}.$
\item ${\Rnp}$ = Space of all $n \times p$ real-valued matrices.
\item $\iMat$ = $p \times p$ identity matrix.
\item $\StiefelS = \{ X \in \Rnp \,:\,X^TX = \iMat\}$, Stiefel Manifold of $p$-frames in $\mathbb{R}_{}^{n}$.
\item $\SpaceM  = \{ X \in \StiefelS : X_{1,j} \geq 0 \;\;\forall\, j=1,2,\cdots,p \}$. 
\item $\SpaceV$ = $O(p)$ = Space of Orthogonal matrices of dimension $p\times p$.
\item $\Measure[]{M}$ = Normalized Haar measure on $\StiefelS$.
\item $\Measure[]{V}$ = Normalized Haar measure on $\SpaceV$.
\item $\Measure[]{D}$ = Lebesgue measure on $\mathbb{R}_{+}^{p}$.

\item $f(\cdot;\cdot)$ = Probability density function.
\item $g(\cdot;\cdot)$ = Unnormalized version of the probability density function.
\item $tr(A)$ = Trace of a square matrix $A$.
\item $etr(A)$ = Exponential of $tr(A)$.
\item ${E}(X)$ = Expectation of the random variable $X$.
\item $\mathbb{I}(\cdot)$ = Indicator function.
\item $\normtwo{\cdot}$ = Matrix operator norm.
\item We use $\bd$ and $D$ interchangeably. $D$ is the diagonal matrix with diagonal $\bd$. We use matrix notation $D$ in the place of $\bd$ wherever needed, and vector $\bd$ otherwise.
\end{itemize}

\section{Properties of the Matrix Langevin distribution and  $\hyp$} \label{Apndx:sec:hypProperties}

We introduce a few lemmas. Readers may skip this section with no loss of understanding of subsequent sections in the paper. 

\begin{lemma}\label{lem:ExpextedValueOfML}
Let $X$ be a random matrix  taking values on the space $\StiefelS$.  If $X\sim \ML(\cdot; M,\bd,V)$, then 
 $E(X)=M D_{\bf h} V^T$. $D_{{ \bf h}}$ is a diagonal matrix with diagonal entrees ${\bf  h}(\bd ):=\left( h_1(\bd), \ldots,  h_p(\bd) \right)$ where
$$ h_j(\bd) :=\frac{ {\frac{\partial }{\partial d_j}\,{}_0F_1 \left(\frac{n}{2}, \frac{D^2}{4} \right)}}{{_0F_1 \left(\frac{n}{2}, \frac{D^2}{4} \right)}} \text{  for } j=1, 2,\cdots, p.$$

\end{lemma}
{\bf Proof of Lemma~\ref{lem:ExpextedValueOfML}.}
Let  $\Gamma_0=\left[ M, \overline{M} \right]$ be a $n\times n$ orthogonal matrix where the columns of the matrix $\overline{M}$ comprise of a orthonormal basis for the orthogonal complement of the  column space of $M$. Consider the random matrix $Y = \Gamma_0^T X V  $. From~\cite{Khatri:1977} (see page 98) we know that 
\begin{eqnarray}\label{eq:ExpectationY}
E(Y)=\left[\begin{array}{c}D_{{\bf h}} \;\;, \;\; \mathbf{0}_{n-p,p} \end{array}\right]^T,
\end{eqnarray}
where 
$D_{{ \bf h}}$ is a diagonal matrix with diagonal entrees ${\bf  h}(\bd):=\left( h_1(\bd), \ldots,  h_p(\bd) \right)$ with
$$ h_j(\bd) :=\frac{ {\frac{\partial }{\partial d_j}\,{}_0F_1 \left(\frac{n}{2}, \frac{D^2}{4} \right)}}{{_0F_1 \left(\frac{n}{2}, \frac{D^2}{4} \right)}} \text{  for } j=1, 2,\cdots, p.$$
 Hence from Equation~\eqref{eq:ExpectationY}, it follows that 
 $$ E(X)=\Gamma_0 \left[\begin{array}{c}D_{{\bf h}} \;\;, \;\; \mathbf{0}_{n-p,p} \end{array}\right]^T V^T = M  D_{{\bf h}} V^T. $$
 \qedwhite
\label{subsec:0F1_prop}

\begin{lemma}
\label{lem:0F1_upper_bound}\citep{Chikuse:2012,Hoff:2009}
For any $p \times p$ diagonal matrix $D$ with positive elements, ${}_0F_1\left( \frac{n}{2}, \frac{D^2}{4}\right) \leq etr(D)$ when $n \geq p$.
\end{lemma}

{\bf Proof of Lemma~\ref{lem:0F1_upper_bound}. }
From Equation~\eqref{m-eq:MLDensity_MDV} in the main article, we have
\begin{eqnarray}
&&\int_{\StiefelS} \MLDensity(X;(M,\bd,V))\;\Measure[1]{X} = 1 \nonumber \\
&\implies& {}_0F_1\left( \frac{n}{2}, \frac{D^2}{4}\right) = \int_{\StiefelS} etr(VDM^T X)\;\Measure[1]{X}.
\label{eq:0F1_integral_form}
\end{eqnarray}
We know that $\MLDensity(X;(M,\bd,V))$ has the {\em{unique modal orientation}} $MV^T$ (page 32 in~\cite{Chikuse:2012}). Hence it follows from Equation~\eqref{eq:0F1_integral_form} that
\begin{eqnarray}
{}_0F_1\left( \frac{n}{2}, \frac{D^2}{4}\right) &\leq& \int_{\StiefelS} etr(VDM^TMV^T)\;\Measure[1]{X} \nonumber \\
&=& etr(D)\int_{\StiefelS} \;\Measure[1]{X} = etr(D),
\label{eq:0F1_upper_bound}
\end{eqnarray}
as  $\Measure[]{X}$ is the normalized Haar measure, i.e. a probability measure on $\StiefelS$.

\qedwhite

\begin{lemma}
\label{lem:diagonal_less_norm}
Let $A$ be a $n \times p$ real matrix with $n\geq p$, and   $A_{j,j}$ be the $(j,j)$-th entry of the matrix $A$ for $j =1, .., p$. Let  $\normtwo{A}$ denote the matrix operator norm ( also known as spectral norm) of the matrix $A$. If  $\normtwo{A} \leq \delta $ for some $\delta>0$ then   $\lvert A_{j,j} \rvert \leq \delta $  for $j =1, .., p$. Also, if   $\normtwo{A} < \delta $ for some $\delta>0$ then   $\lvert A_{j,j} \rvert < \delta $  for $j =1, .., p$.
\end{lemma}

{\bf Proof of Lemma~\ref{lem:diagonal_less_norm}.}

From the assumptions of Lemma~\ref{lem:diagonal_less_norm} along with the definition of the spectral norm, it follows that  $l^T A^T A\; l \leq \delta^2 $ for all $l\in \mathbb{R}^p$ with $l^T l=1$. In particular, $e_j^T A^T A\; e_j \leq \delta^2 $ where $e_j\in \mathbb{R}^p$ such that its  $j$-{th} entry equals $1$ while rest  of its entries are $0$. Hence we have $\sum_{k=1}^{n}A_{k,j}^2 \leq \delta^2 $ implying the fact that $|A_{j,j}| \leq \delta$. The assertion with strict inequality can also be shown in a similar fashion.

\qedwhite 
\begin{lemma}
\label{lem:F_lowerbound} 
Let $D$ be a $p \times p$ diagonal matrix with positive diagonal elements $\bd = (d_1, d_2, \cdots, d_p)$. Then for any $\delta>0$ and $n \geq p$, there exists a positive constant, $\ConstOne_{n,p,\delta}$, depending on $n, p$ and $\delta$, such that 
$$_0F_1\left(\frac{n}{2}, \frac{D^2}{4}\right) > \ConstOne_{n,p,\delta} \;etr\left((1-\delta)D\right).$$
\end{lemma}

{\bf Proof of Lemma~\ref{lem:F_lowerbound}.}

Note that $D$ is a $p\times p$ diagonal matrix with positive diagonal elements   $d_1, .., d_p$.  For the case $n \geq p$, define 

\begin{equation}
\widetilde{M} = \left[
\begin{array}{c}
\iMat\;\;\;\; \\
\mathbf{0}_{n-p,p} 
\end{array}
\right],
\widetilde{V} = \iMat
\;\;\mbox{and}\;\;
I^{\star}:=
\left[
\begin{array}{c}
\iMat \;\;\;\; \\
\mathbf{0}_{n-p,p} 
\end{array}
\right],
\label{eq:defn_M_tilde_V_tilde}
\end{equation}

where $\iMat$ denotes the $p\times p$ identity matrix and $\mathbf{0}_{n-p,p}$ represents the zero matrix of dimension $(n-p) \times p$. For arbitrary given positive constant $\delta>0$, consider
$$ B_{\delta}:=\left\{X \in \mathcal{V}_{n,p }, \text{ such that }   \normtwo{X- I^{\star}}  <\delta\right\},$$

where $\normtwo{\cdot}$ denotes the spectral norm of a matrix. Let $\HaarMu$ denotes the normalized Haar measure on the $\mathcal{V}_{n,p}$. Clearly,  $0<\HaarMu \left(B_{\delta}\right)<\infty $, as $B_{\delta}$ is a non-empty open subset of $\mathcal{V}_{n,p }$. 
Now from Equation~\eqref{m-eq:MLDensity_MDV} we have, 
\begin{eqnarray}
_0F_1\left(\frac{n}{2}, \frac{D^2}{4}\right)& = & \int_{\mathcal{V}_{n,p}} etr\left(\widetilde{V}D\widetilde{M}^T X\right)d\HaarMu(X).\nonumber\\
& \geq  & \int_{B_{\delta}} etr\left(\widetilde{V}D\widetilde{M}^T X\right)d\HaarMu(X).
\label{eq:ZeroFone_lower_1}
\end{eqnarray}
 

Using Lemma~\ref{lem:diagonal_less_norm} we know that $X_{{j,j}}>(1-\delta)$ for $j = 1, 2, \ldots, p$ where $X \in B_{\delta}$. Note that $X_{{j,j}}$ denotes the $(j,j)$-th entry of the matrix $X$. Hence from Equation~\eqref{eq:defn_M_tilde_V_tilde} and~\eqref{eq:ZeroFone_lower_1} it follows that, 
\begin{eqnarray}
\label{eq:ZeroFone_lower_3}
_0F_1\left(\frac{n}{2}, \frac{D^2}{4}\right) 
& \geq  & \int_{B_{\delta}} \exp\left({\sum_{j=1}^p {X}_{j,j}\,d_j}\right) \,d\HaarMu(X), \nonumber \\
&>&  \HaarMu(B_{\delta}) \; etr\left((1-\delta)D\right),
\end{eqnarray}

where the last inequality uses the fact that $d_j>0$ for all $j =1, \ldots p.$ Finally we denote  $ \ConstOne_{n,p,\delta}:= \HaarMu(B_{\delta})>0$ as it depends on  $n,p$ along with $\delta$, to conclude that 

$$_0F_1\left(\frac{n}{2}, \frac{D^2}{4}\right) > \ConstOne_{n,p,\delta} \;etr\left((1-\delta)D\right).$$

\qedwhite

\begin{lemma}
For any $p \times p$ diagonal matrix $D$ with positive elements $\bd \in \SpaceD$, the hypergeometric function of matrix argument denoted by $\hyp$ is log-convex with respect to $\bd$ where $n \geq p$.
 
\label{lem:0F1_log_convex}
\end{lemma}

{\bf Proof of Lemma~\ref{lem:0F1_log_convex}.}\\
From Equation~\eqref{m-eq:MLDensity_MDV} in the main article, we have
\begin{equation}
\hyp = \int_{\StiefelS} etr(VDM^TX)\,\Measure[1]{X}, 
\label{eq:0F1_int_form}
\end{equation}
for arbitrary $M \in \SpaceM$ and $V \in \StiefelS$ where $n \geq p$. Without loss of generality, we can take $M = \widetilde{M} = \left[ \begin{array}{l}
\iMat \\
\bm{0}_{(n-p),p} \end{array} \right]$ 
and $V = \iMat$.

Let $D_1$ and $D_2$ be two $p \times p$ diagonal matrices with positive diagonal entries $\bd_1$ and $\bd_2$, respectively, and $\bd_1 \neq \bd_2$. From Equation~\eqref{eq:0F1_int_form}, we have
\begin{eqnarray}
{}_0F_1\left(\frac{n}{2},\frac{D_1^2}{4}\right) =  \int_{\StiefelS} etr(D_1\widetilde{M}^TX)\,\Measure[1]{X} \nonumber \\
{}_0F_1\left(\frac{n}{2},\frac{D_2^2}{4}\right) = \int_{\StiefelS} etr(D_2\widetilde{M}^TX)\,\Measure[1]{X}.
\label{eq:0F1_D1_D2}
\end{eqnarray}
Let $\lambda \in (0,1)$ be any real number, then
\begin{eqnarray}
&&{}_0F_1\left(\frac{n}{2},\frac{{\left(\lambda\,D_1 + (1-\lambda)\,D_2\right)}^2}{4}\right) \nonumber \\
&=&\int_{\StiefelS} etr((\lambda\,D_1 + (1-\lambda)\,D_2)\tilde{M}^TX)\,[dX] \nonumber \\
&=&\int_{\StiefelS} {\left(etr(D_1\tilde{M}^TX)\right)}^\lambda {\left(etr(D_2\tilde{M}^TX)\right)}^{1-\lambda}\,\Measure[1]{X} \nonumber \\
&<& {\left(\int_{\StiefelS}\,etr(D_1\tilde{M}^TX)\,\Measure[1]{X}\right)}^\lambda {\left(\int_{\StiefelS}\,etr(D_2\tilde{M}^TX)\,\Measure[1]{X}\right)}^{1-\lambda} \nonumber \\
&=&{\left({}_0F_1\left(\frac{n}{2},\frac{D_1^2}{4}\right)\right)}^\lambda\,{\left({}_0F_1\left(\frac{n}{2},\frac{D_2^2}{4}\right)\right)}^{1-\lambda}.
\label{eq:convexity_0F1}
\end{eqnarray}
where the inequality is due to  H\"older~\citep{Hardy:1952,Book:Billingsley:1995} and we have strict inequality as  $\bd_1 \neq \bd_2$. 
Therefore from Equation~\eqref{eq:convexity_0F1}, it follows that
$\hyp$ is a log-convex function of the diagonal entries $\bd$ of the matrix $D$. Note that, the properties of the exponential family of distributions have played a crucial role in establishing the result. 

\qedwhite


\begin{lemma}
For any $p \times p$ ($p\geq 2$) diagonal matrix $D$ with positive elements $\bd \in \SpaceD$, 
$$ 0 < \; \partialhyp{i} < \hyp,$$ 
for $i=1, 2, \cdots, p$, where $n \geq p$.
\label{lem:0F1_prime}
\end{lemma}

{\bf Proof of Lemma~\ref{lem:0F1_prime}.}

\paragraph{Right hand side inequality} Proceeding along similar lines as Lemma~\ref{lem:0F1_log_convex} we have
\begin{equation}
\hyp = \int_{\StiefelS} etr(D\widetilde{M}^TX)\,\Measure[1]{X}, \;\;\mbox{where $\widetilde{M} = \left[ \begin{array}{l}
\iMat \\
\bm{0}_{(n-p),p} \end{array} \right]$}. 
\label{eq:0F1_int_form_tilde_M}
\end{equation}
From Equation~\eqref{eq:0F1_int_form_tilde_M}, we have 
\begin{eqnarray}
\hyp &=& \int_{\StiefelS} \exp\left(\sum_{j=1}^p d_j\,X_{j,j}\right)\,\Measure[1]{X}
\label{eq:0F1_int_form_X_i_0}
\end{eqnarray}
Consider the set $\mathcal{V}_{0} := \left\{ X \in \StiefelS : X_{i,i}=1 \right\}.$ 
Note that $\mathcal{V}_{0}$ is  isomorphic  to the lower dimensional Stiefel manifold, $\mathcal{V}_{n,p-1}$. $\mathcal{V}_{0}$, being a lower dimensional subspace of $\StiefelS$, has measure zero i.e. $\int_{\StiefelS }\IndVzero{} \Measure[1]{X}=0$, where $\IndVzero{}$ is the indicator function for $X$ to be in the set $\mathcal{V}_0$. 
From Equation~\eqref{eq:0F1_int_form_X_i_0}, we have 
\begin{eqnarray}
\hyp &=& \int_{\StiefelS} \exp\left(\sum_{j=1}^p d_j\,X_{j,j}\right)\,\IndVzero{c}\;\Measure[1]{X},
\label{eq:0F1_new_form}
\end{eqnarray}
where $\mathcal{V}^c_0$ is the complement of $\mathcal{V}_0$. Hence,


\begin{eqnarray}
\partialhyp{i} &=& \int_{\StiefelS} X_{i,i}\,\IndVzero{c}\,\exp\left(\sum_{j=1}^p d_j\,X_{j,j}\right)\,\Measure[1]{X}.\nonumber\\ 
\label{eq:0F1_int_form_X_i_1}
\end{eqnarray}
Observe that $\normtwo{X}=1$ on $\StiefelS$. Hence from Lemma~\ref{lem:diagonal_less_norm}  we have $\lvert  X_{i,i} \rvert  \leq 1$. Also,  $ X_{i,i} \neq 1$ when $X \in \mathcal{V}^{c}_{0}$. As a result, we conclude that $X_{i,i} < 1$ on ${\StiefelS \cap \mathcal{V}^{c}_{0} }$. Consequently, it follows from Equations~\eqref{eq:0F1_new_form} and~\eqref{eq:0F1_int_form_X_i_1} that, 
\begin{eqnarray}
\partialhyp{i}&\stackrel{}{<} & \int_{\StiefelS} \exp\left(\sum_{j=1}^p d_j\,X_{j,j}\right)\,\IndVzero{c}\,\Measure[1]{X} \nonumber \\
& = & \hyp. 
\label{eq:0F1_int_form_X_i_2}
\end{eqnarray}


\paragraph{Left hand side inequality} 
Consider  $\StiefelS^{i,+}:=\left\{ X \in \StiefelS : X_{i,i}>0 \right\}$, $\StiefelS^{i,-}:=\left\{ X \in \StiefelS : X_{i,i}<0 \right\}$ and $\StiefelS^{i,0}:=\left\{ X \in \StiefelS : X_{i,i}=0 \right\}$. 
Clearly, $\StiefelS^{i,+} , \StiefelS^{i,0}$ and $\StiefelS^{i,-}$ forms a partition of $\StiefelS$. Hence from Equation~\eqref{eq:0F1_int_form_X_i_0} we have, 
\begin{eqnarray}
 & & \partialhyp{i} \nonumber\\
 &=& \int_{\StiefelS^{i,+}} X_{i,i} \exp\left(\sum_{j=1}^p d_j\,X_{j,j}\right)\,\Measure[1]{X}+ \int_{\StiefelS^{i,0}}  X_{i,i}\exp\left(\sum_{j=1}^p d_j\,X_{j,j}\right)\,\Measure[1]{X} \nonumber \\
 &&\hspace{2in}+\int_{\StiefelS^{i,+}} X_{i,i}\exp\left(\sum_{j=1}^p d_j\,X_{j,j}\right)\,\Measure[1]{X}\nonumber\\
 &=& \int_{\StiefelS^{i,+}} X_{i,i} \exp\left(\sum_{j=1}^p d_j\,X_{j,j}\right)\,\Measure[1]{X}+ \int_{\StiefelS^{i,-}}  X_{i,i}\exp\left(\sum_{j=1}^p d_j\,X_{j,j}\right)\,\Measure[1]{X}.\nonumber \\
\label{eq:0F1_int_form_X_i_5}
\end{eqnarray}

Let $\Gamma_{}$ be the  $n \times n $ diagonal matrix such that $\Gamma_{j,j}=1$ for $j= 1, \ldots , n, j \neq i$ and $\Gamma_{i,i}=-1$.  $\Gamma_{}$ is an orthogonal matrix as $\Gamma^T\Gamma = \bf{I}_n$. It is easy to show that $\StiefelS^{i,+}= \left\{ \Gamma\, X : X \in \StiefelS^{i,-} \right\}$.

Consider the change of variable $Y:= \Gamma X$.  Using standard algebra we can show that $X_{i,i}=-Y_{i,i}$ and $X_{j,j}=Y_{j,j}$ for $j =1, \ldots p, j \neq i $. As the normalized Haar  measure on $\StiefelS$ is invariant under orthogonal transformation from~\cite{Chikuse:2012}, we get that  
\begin{eqnarray}
& & \int_{\StiefelS^{i,-}}  X_{i,i}\exp\left(\sum_{j=1}^p d_j\,X_{j,j}\right)\,\Measure[1]{X}\nonumber\\
& =&  - \int_{\StiefelS^{i,+}}  Y_{i,i}\exp\left(-d_i\,Y_{i,i}+\sum_{j=1, j \neq i}^p d_j\,Y_{j,j}\right)\,\Measure[1]{Y} \nonumber\\
 &=& - \int_{\StiefelS^{i,+}}  X_{i,i}\exp\left(-d_i\,X_{i,i}+\sum_{j=1, j \neq i}^p d_j\,X_{j,j}\right)\,\Measure[1]{X}.
\label{eq:0F1_int_form_X_i_5_1}
\end{eqnarray}

From Equations~\eqref{eq:0F1_int_form_X_i_5}  and~\eqref{eq:0F1_int_form_X_i_5_1} we have,

\begin{eqnarray}
& & \partialhyp{i} \nonumber\\
&=& \int_{\StiefelS^{i,+}} X_{i,i} \, \exp\left(\sum_{j=1, j \neq i}^p d_j\,X_{j,j}\right) \Bigg(\exp\left( d_i\,X_{i,i}\right)- \exp\left(-d_i\,X_{i,i}\right)\Bigg)\,\Measure[1]{X}\nonumber\\
&=& \int_{\StiefelS^{i,+}} X_{i,i}  \exp\left( \sum_{j=1, j \neq i}^p d_j\,X_{j,j}\right) 2\, \sinh\left( d_i\,X_{i,i}\right) \,\Measure[1]{X}\nonumber\\
\label{eq:0F1_int_form_X_i_6}
\end{eqnarray}

where $\sinh$ is the hyperbolic sin function. Note that $\sinh\left( d_i\,X_{i,i}\right)>0$ as  $ d_i>0$ and $X_{i,i}>0$ on $\StiefelS^{i,+}$ . Hence from Equation~\eqref{eq:0F1_int_form_X_i_6} it follows that, 
\begin{eqnarray}
\partialhyp{i}  &\stackrel{}{>}& 0.
\label{eq:derivative_postivity}
\end{eqnarray}

From Equations~\eqref{eq:0F1_int_form_X_i_2} and~\eqref{eq:derivative_postivity}, we have the result.

\qedwhite
\newline
\begin{lemma}
\label{lem:prior_DY1}
Let $\priorXzero \in \Rnp$ and  $D$ be  a diagonal matrix with positive diagonal entries. If $\normtwo{\priorXzero} < 1$, then for arbitrary  $M\in \StiefelS, V \in \SpaceV$,  
\begin{eqnarray}
\frac{etr\left( VDM^T\priorXzero\right)}{ _0 F_1 (\frac{n}{2}, \frac{D^2}{4})} < \frac{etr(-\epsilon_0\,D)}{\ConstOne_{n,p,\epsilon_0}},
\end{eqnarray}
where $\epsilon_0=\frac{1}{2}\left(1-\normtwo{\priorXzero}\right)$ and $\ConstOne_{n,p,\epsilon_0} >0 $ is a constant depending on $n,p$ and $\epsilon_0$.
\end{lemma}

{\bf Proof of Lemma~\ref{lem:prior_DY1}.}

Note that $0<\epsilon_0<\frac{1}{2}$, as $\normtwo{\priorXzero} < 1$. Assume $Y_0= M^T \priorXzero V \in \Rpp$. For arbitrary $l\in \mathbb{R}^p$ with $\vecnorm{l}=1$, we have
\begin{eqnarray}
l^T  Y_0^T Y_0 l &=&  (V\,l)^T  \priorXzero^T \priorXzero (V\,l) -l^T  V^T \priorXzero^T ({\bf{I}}_{n}-M M^T) \priorXzero V l \nonumber \\
&\leq&( 1-2 \epsilon_0)^2.
\label{eq:YnormIsLess}
\end{eqnarray}
The last inequality follows as $\normtwo{\priorXzero} = 1-2\epsilon_0$ and $({\bf{I}}_{n}-MM^T) $ is a non-negative definite matrix. From Equation~\eqref{eq:YnormIsLess} it follows that $\normtwo{Y_0} \leq 1-2\epsilon_0$. Hence, applying Lemma~\ref{lem:diagonal_less_norm} we obtain that $\lvert{Y_0}_{j,j}\rvert < 1-2\epsilon_0$ for $j=1,\cdots, p$, where $Y_{0,j}$ is the $j$-{th} diagonal element of the matrix $Y_0$. Now applying Lemma~\ref{lem:F_lowerbound} we have,
\begin{eqnarray}
\frac{etr\left( VDM^T \priorXzero\right)}{ _0 F_1 (\frac{n}{2}, \frac{D^2}{4})} &<& \frac{etr(D Y_0-(1-\epsilon_0)D)}{\ConstOne_{n,p,\epsilon_0}}< \frac{etr(-\epsilon_0\,D)}{\ConstOne_{n,p,\epsilon_0}}. \nonumber
\end{eqnarray}
\qedwhite
 \newline

{
\begin{lemma}\label{lem:LowerBoundOfHyp}
Let $R$ be a $p\times p$ symmetric positive definite matrix. Then for $a\geq p/2$,
\begin{eqnarray}
{}_0F_1\left(a,R\right)\geq  \Gamma(a) ({tr(R)})^{\frac{1-a}{2}}I_{a-1}(\sqrt{4\;tr(R) }),
\end{eqnarray}
$tr(R)$ denotes the trace of the matrix $R$.
\end{lemma}
{\bf Proof of Lemma~\ref{lem:LowerBoundOfHyp}.}
Let $\mathcal{D}_{k,p}$ denotes the set of all possible partitions of the integer $k$ into no more than $p$ parts, i.e. $$\mathcal{D}_{k,p}=\left\{(k_1,\ldots, k_p): k_1, \ldots k_p \in \mathbb{Z}, k_1\geq \ldots k_p\geq 0, k_1+ \ldots +k_p=k  \right\},$$
where $\mathbb{Z}$ denotes the set of non-negative integers. For a vector $\boldsymbol{\kappa}=(k_1,\ldots, k_p)\in \mathcal{D}_{k,p}$, we denote the quantity  $\prod_{j=1}^{p} \frac{\Gamma(a-(j-1)/2 +k_j)}{\Gamma(a-(j-1)/2)}$  by the notation $(a)_{\boldsymbol{\kappa}}$. Then from the \cite{Richards:2011} we get the representation  
\begin{eqnarray}
{}_0F_1\left(a,R\right)=\sum_{k=0}^{\infty}\sum_{\boldsymbol{\kappa}\in \mathcal{D}_{k,p}}\frac{ C_{\boldsymbol{\kappa}}(R)}{(a)_{\boldsymbol{\kappa}} k!},
\end{eqnarray}
where $C_{\boldsymbol{\kappa}}(R)$ is the Zonal polynomial of the matrix argument $R$ corresponding to the vector $\boldsymbol{\kappa}\in\mathcal{D}_{k,p}$. More details about the Zonal polynomials can be found in \cite{Muirhead:2009:aspects}, \cite{Gross:1987}. 

Note that, for ${\boldsymbol{\kappa}=(k_1, \ldots, k_p)\in \mathcal{D}_{k,p}}$, and $a\geq \frac{p}{2}$,
\begin{eqnarray}
(a)_{\boldsymbol{\kappa}}&=&\prod_{j=1}^{p} \frac{\Gamma(a-(j-1)/2 +k_j)}{\Gamma(a-(j-1)/2)}\nonumber\\
&\leq &  \frac{\Gamma(a +k)}{\Gamma(a)} \prod_{j=2}^{p} \frac{\Gamma(a-(j-1)/2 )}{\Gamma(a-(j-1)/2)}\nonumber\\
&= &  \frac{\Gamma(a +k)}{\Gamma(a)}.
\end{eqnarray} 
As a result, for $a\geq \frac{p}{2}\geq 1$, we get that 
\begin{eqnarray}
{}_0F_1\left(a,R\right)
&\geq& \sum_{k=0}^{\infty}\frac{\Gamma(a)}{k!\Gamma(a +k)}\sum_{\boldsymbol{\kappa}\in \mathcal{D}_{k,p}}{ C_{\boldsymbol{\kappa}}(R)}\nonumber\\
&\stackrel{(**)}{=}& \sum_{k=0}^{\infty}\frac{\Gamma(a)}{k!\Gamma(a +k)} tr(R)^k\nonumber\\
&\stackrel{ }{=}& \Gamma(a) ({tr(R)})^{\frac{1-a}{2}}I_{a-1}(\sqrt{4\;tr(R) }),
\end{eqnarray}
where the equality in $(**)$ follows from \cite{Gross:1987} (See Equation(5) in \cite{Gross:1987} ), while the last equality follows from the definition of $I_{a-1}(\cdot)$, the modified Bessel function of the first kind. We would like to point out that the result is motivated by a lower-bound developed in \cite{sengupta2013two}.
\qedwhite

\begin{lemma}\label{lem:BesselIProperty}
Let $\nu\geq \frac{1}{2}$ then for $M>0$, 
\begin{eqnarray}
 I_{\nu}(x) \geq {\frac{e^{x}}{\sqrt{x}}}  \left[\sqrt{M}e^{-M}I_{\nu}(M)\right],
\end{eqnarray}
for all $x>M$.
\end{lemma}

{\bf Proof of Lemma~\ref{lem:BesselIProperty}}

First we will show that the function $x\mapsto x^{\frac{1}{2}}e^{-x}I_{\nu}(x)$ is a non decreasing function for $\nu\geq \frac{1}{2}$ and $x>0$. Consider that

%

\begin{eqnarray}\label{eq:IbesselIneq1}
& &\frac{\partial}{\partial x} (x^{\frac{1}{2}}e^{-x}I_{\nu}(x))\nonumber\\
 &=& \frac{1}{2\sqrt{x}}e^{-x}I_{\nu}(x)-(x^{\frac{1}{2}}e^{-x}I_{\nu}(x))+x^{\frac{1}{2}}e^{-x} (-\frac{\nu}{ x}I_{\nu}(x)+ I_{\nu-1}(x))\nonumber\\
 &=& {\sqrt{x}}e^{-x} \left(  (\frac{1}{2}-\nu)\frac{I_{\nu}(x)}{x}+ I_{\nu-1}(x)-I_{\nu}(x)\right)\nonumber\\
 &=& {\sqrt{x}}e^{-x} I_{\nu}(x)\left(  \frac{0.5-\nu}{x}-1+ \frac{I_{\nu-1}(x)}{I_{\nu}(x)}\right).
\end{eqnarray}

From \cite{Segura:2011}, we get that $\frac{I_{\nu}(x)}{I_{\nu-1}(x)}\leq \frac{x}{(\nu-0.5)+\sqrt{(\nu-0.5)^2+x^2}}$ for $\nu\geq \frac{1}{2}$ and $x>0$. Hence, from \eqref{eq:IbesselIneq1}, it follows that

\begin{eqnarray}
& &\frac{\partial}{\partial x} (x^{\frac{1}{2}}e^{-x}I_{\nu}(x))\nonumber\\
  &\geq & {\sqrt{x}}e^{-x} I_{\nu}(x)\left(  \frac{0.5-\nu}{x}-1+ \frac{(\nu-0.5)+\sqrt{(\nu-0.5)^2+x^2}}{x}\right)\nonumber\\
   &= & {\sqrt{x}}e^{-x} I_{\nu}(x)\left(   \frac{\sqrt{(\nu-0.5)^2+x^2}}{x}-1\right)\nonumber\\
    &= &   \frac{{\sqrt{x}}e^{-x} I_{\nu}(x) (\nu-0.5)^2}{x\left(x+\sqrt{(\nu-0.5)^2+x^2}\right)}>0.\nonumber
\end{eqnarray}

 As a result, the function $x\mapsto (x^{\frac{1}{2}}e^{-x}I_{\nu}(x))$ is a non-decreasing function for $\nu\geq \frac{1}{2}$.  Hence, for $M>0$ we have 
\begin{eqnarray}
& & x^{\frac{1}{2}}e^{-x}I_{\nu}(x) \geq M^{\frac{1}{2}}e^{-M}I_{\nu}(M)\nonumber\\
& \implies&  I_{\nu}(x) \geq {\frac{e^{x}}{\sqrt{x}}}  \left[\sqrt{M}e^{-M}I_{\nu}(M)\right],
\end{eqnarray}
when $x>M$.
\qedwhite

\begin{lemma}\label{lem:g1UpperBound}
Let $n > p \geq 2$ and $M>0$, 
for all $d_1>M$,

\begin{eqnarray}
g_1(d_1)&\leq & K_{n,p,M} \;d_1^{\nu(n-1)/2} exp(\;-\nu(1-\eta_1)\,d_1),\nonumber\\
\end{eqnarray}
where $$g_1(d_1)=\frac{exp(\nu\;\eta_1\,d_1)}{\hyp^{\nu}} \text{ and } K^{\dagger}_{n,p,M}= \left[ \frac{({p /4})^{\frac{n/2-1}{2}})}{\Gamma(n/2)\left\{  \sqrt{M}e^{-M}I_{n/2-1}(M)\right\} } \right]^{\nu}.$$
\end{lemma}

{\bf Proof of Lemma~\ref{lem:g1UpperBound}}
Let $a=\frac{n}{2}$. Note that $d_1 \geq d_2\geq \ldots \geq d_p$ are the diagonal elements of the diagonal matrix $D$. 
From Lemma~\eqref{lem:LowerBoundOfHyp}, we get that

\begin{eqnarray}
{}_0F_1\left(a,\frac{D^2}{4}\right)&\geq&  \Gamma(a) \left(\frac{4}{tr(D^2)}\right)^{\frac{a-1}{2}}I_{a-1}(\sqrt{tr(D^2) })\nonumber\\
&\geq&  \Gamma(a) \left(\frac{4}{p d_1^2}\right)^{\frac{a-1}{2}}\;I_{a-1}(d_1),\nonumber\\
\end{eqnarray}

As a result, 
\begin{eqnarray}\label{eq:g1Upperbound0}
g_1(d_1)&\leq & \left[ \frac{({p d_1^2/4})^{\frac{a-1}{2}} exp(\;\eta_1\,d_1)}{\Gamma(a)I_{a-1}(d_1) } \right]^{\nu}.\nonumber\\
\end{eqnarray}

With the help of Lemma~\eqref{lem:BesselIProperty}, from Equation~\eqref{eq:g1Upperbound0}, we get that
\begin{eqnarray}\label{eq:g1upperbound}
g_1(d_1)&\leq & \left[ \frac{({p d_1^2/4})^{\frac{a-1}{2}}) exp(\;\eta_1\,d_1)}{\Gamma(a)\left\{ {\frac{e^{d_1}}{\sqrt{d_1}}}  \left[\sqrt{M}e^{-M}I_{a-1}(M)\right]\right\} } \right]^{\nu}\nonumber\\
& =& \left[ \frac{({p /4})^{\frac{a-1}{2}})d_1^{a-0.5} exp(\;-(1-\eta_1)\,d_1)}{\Gamma(a)\left\{  \left[\sqrt{M}e^{-M}I_{a-1}(M)\right]\right\} } \right]^{\nu}\nonumber\\
& =& \left[ \frac{({p /4})^{\frac{a-1}{2}})}{\Gamma(a)\left\{  \sqrt{M}e^{-M}I_{a-1}(M)\right\} } \right]^{\nu} d_1^{\nu(a-0.5)} exp(\;-\nu(1-\eta_1)\,d_1).\nonumber\\
\end{eqnarray}

Note that $lim_{M\rightarrow \infty} \sqrt{M} e^{-M}I_{a-1}(M)=\frac{1}{\sqrt{2\pi}}$ for all $a\geq \frac{3}{2}.$ The upper bound is nontrivial in the sense that for larger values of $M$, the constant part involved in the inequality \eqref{eq:g1upperbound} does not approach infinity.
} 
 
All the above lemmas are used for the theoretical development of the Bayesian analysis with $\ML$ distributions.



\section{The matrix Langevin distribution on the Stiefel manifold}
\label{sec:stiefel_distr}

The Stiefel manifold, $\StiefelS$, is the space of all $p$ ordered orthonormal vectors (also known as $p$-frames) in $\mathbb{R}^n$ \citep{Mardia:2009, Absil:2009, Chikuse:2012, Edelman:1998, Downs:1972} and is defined as 
\begin{equation*}
\StiefelS = \{ X \in \Rnp \,:\,X^TX = \iMat,\;  p\leq n\}, 
\label{eq:manifold_stiefel}
\end{equation*}
where $\Rnp$ is the space of all $n \times p$, $p\leq n$ real-valued matrices, and $\iMat$ is the $p \times p$ identity matrix. $\StiefelS$ is a compact Riemannian manifold of dimension $np - p(p + 1)/2$ \citep{Chikuse:2012}. A topology on $\StiefelS$ can be induced from the topology  on $\Rnp$ as $\StiefelS$ is a sub-manifold of  $\Rnp$ ~\citep{Absil:2009,Edelman:1998}.
For $p=n$, $\StiefelS$ becomes identical to $O(n)$, the orthogonal group consisting of all orthogonal $n \times n$ real-valued matrices, with the group operation being matrix multiplication.
Being a compact unimodular group, $O(n)$ has a unique Haar measure that corresponds to a uniform probability measure on $O(n)$ \cite{Chikuse:2012}. Also, through obvious mappings, the Haar measure on $O(n)$ induces a normalized Haar measure on the compact manifolds $\StiefelS$. The normalized Haar measures on $O(n)$ and $\StiefelS$ are invariant under orthogonal transformations \citep{Chikuse:2012}.
Detailed construction of the Haar measure on $\StiefelS$ and its properties are described in~\cite{Muirhead:2009, Chikuse:2012}. Notation wise, we will use $\Measure[]{M}$ and $\Measure[]{V}$ to denote the normalized Haar measures  on $\StiefelS$ and $\SpaceV$, respectively.

The matrix Langevin distribution ($\ML$-distribution) is a widely used probability distribution on $\StiefelS$~\citep{Mardia:2009,Chikuse:2012,Lin:2017}. This distribution is also known as Von Mises-Fisher matrix distribution~\citep{Khatri:1977}. 
As defined in~\cite{Chikuse:2012}, the probability density function of the matrix Langevin distribution (with respect to the normalized Haar measure $\Measure[]{M}$ on $\StiefelS$) parametrized by $F \in \Rnp$, is 
\begin{equation}
\MLDensity (X\,;\,F) = \frac{etr(F^T X)}{{}_0F_1\left( \frac{n}{2}, \frac{F^TF}{4}\right)}, 
\label{eq:MLDensity_F}
\end{equation}
where $etr(\cdot) = \exp(trace(\cdot))$  and the normalizing constant, ${}_0F_1({n}/{2}, {F^TF}/{4})$, is the hypergeometric function of order $n/2$ with the matrix argument $F^TF/4$~\citep{Herz:1955,James:1964,Muirhead:1975,Gupta:1985,Gross:1987,Gross:1989,Butler:2003,Koev:2006,Chikuse:2012}. In this article, we consider a different parametrization of the parameter matrix $F$ in terms of its singular value decomposition (SVD). In particular, we subscribe to the specific form  of unique SVD defined in \cite{Chikuse:2012} (Equation 1.5.8 in~\cite{Chikuse:2012}), 
%
%
$$F = MDV^T$$
 where $M \in \SpaceM$, $V \in \SpaceV$, and $D$ is the diagonal matrix with diagonal entries $\bd = (d_1, d_2, \cdots, d_p) \in \SpaceD$. Here $\SpaceM  = \{ X \in \StiefelS : X_{1,j} \geq 0 \;\;\forall\, j=1,2,\cdots,p \}$  and $\mathcal{S}_p= \left\{\left(d_1, \ldots, d_p\right) \in \mathbb{R}_{+}^{p} : 0< d_p< \cdots < d_1 <\infty \right\}.$  Henceforth, we shall use the phrase ``unique SVD" to refer to this specific form of SVD.
\noindent
~\cite{Khatri:1977} (page 96) shows that the function  ${}_0F_1({n}/{2}, {F^TF}/{4})$ depends only on the eigenvalues of the matrix $F^TF$, i.e., 
\begin{equation*}
{}_0F_1\left(\frac{n}{2},\frac{F^TF}{4}\right)={}_0F_1\left(\frac{n}{2},\frac{D^2}{4}\right).
\end{equation*}

As a result, we reparametrize the $\ML$ density as 
\begin{equation*}
\MLDensity (X;(M,\bd,V)) = \frac{etr(VDM^T X)}{{}_0F_1(\frac{n}{2}, \frac{D^2}{4})}\;\mathbb{I}(M \in \SpaceM,\bd \in \SpaceD, V \in \SpaceV).
\end{equation*}
This parametrization ensures identifiability of all the parameters $M, \bd$ and $V$. With regard to interpretation, the mode of the distribution is $M V^T$ and $\bd$ represents the concentration parameter~\citep{Chikuse:2003}. For notational convenience we omit the indicator function and write the $\ML$ density as
\begin{equation}
\MLDensity (X;(M,\bd,V)) = \frac{etr(VDM^T X)}{{}_0F_1( \frac{n}{2}, \frac{D^2}{4})}.
\label{eq:MLDensity_MDV}
\end{equation}
where it is understood that $M \in \SpaceM,\bd \in \SpaceD, V \in \SpaceV$. The parametrization with $M, \bd$ and $V$ enables us to represent the intractable hypergeometric function of a matrix argument as a function of vector $\bd$, the diagonal entries of $D$, paving a path for an efficient posterior inference procedure.

We note in passing that an alternative parametrization through polar decomposition with $F=MK$~\citep{Mardia:2009} may pose computational challenges since the elliptical part $K$ lies on a positive semi-definite cone and inference on positive semi-definite cone is not straightforward~\citep{Hill:1987,Bhatia:2007,Schwartzman:2006}. 


\section{Conjugate Prior for the $\ML$-Distribution}
\label{sec:prior_construct}

In the context of the exponential family of distributions, \cite{Diaconis:Ylvisaker:1979} (DY) provides a standard procedure to obtain a class of conjugate priors when the distribution is represented through its natural parametrization~\citep{Casella:2002}. Unfortunately, for the $\ML$ distribution, the DY theorem can not be applied directly, as demonstrated next. We therefore develop, in Section~\ref{subsec:prior_novel}, two novel classes of priors and present a detailed investigation of their properties.

\subsection{Inapplicability of DY theorem for construction of priors for the $\ML$-distribution}
\label{subsec:DY_inapplicable}


In order to present the arguments in this section, we introduce notations  $P_{\theta}$, $x_{_A}$, $\Measure[]{B}$, and ${\Measure[]{B}}_A$,  that are directly drawn from ~\cite{Diaconis:Ylvisaker:1979}. In brief, $P_{\theta}$ denotes the probability measure that is absolutely continuous with respect to an appropriate $\sigma$-finite measure $\Measure[]{B}$ on a convex subset of the Euclidean space, $\mathbb{R}^d$.
 In the case of the $\ML$ distribution, $\Measure[]{B}$ is the Haar measure defined on the Stiefel manifold. The symbol $\mathcal{X}$  denotes the interior of the support of the measure $\Measure[]{B}$. As shown in~\cite{Hornik:2013} $\mathcal{X}:=\left\{ X: \normtwo{X}<1\right\}$ for the case of the $\ML$ distribution. According to the assumptions of DY $\int_{\mathcal{X}}dP_{\theta}(X)=1$ (see paragraph after equation (2.1), page 271 in~\cite{Diaconis:Ylvisaker:1979}). In the current context,  $P_{\theta}$ is the probability measure associated with the $\ML$ distribution. Therefore, 
$$ \int_{\mathcal{X}}dP_{\theta}(X) =\int_{\mathcal{X}}f_{\ML}\left( X\right)\Measure[]{B}(dX)=0,$$
which violates the required assumption mentioned above.
Secondly, in the proof of Theorem $1$ in~\cite{Diaconis:Ylvisaker:1979} DY construct a probability measure restricted to a measurable set $A$ as follows.
$${\Measure[]{B}}_A(B)=\frac{{\Measure[]{B}}(A\cap B)}{{\Measure[]{B}}(A)}, \mbox{ where }{\Measure[]{B}}(A)>0. $$
Considering the notation $x_{_A} = \int Z \;{\Measure[]{B}}_A(dZ)$ for any measurable set $A$,
the proof of Theorem 1 in~\cite{Diaconis:Ylvisaker:1979} relies on the existence of a sequence of measurable sets $\left\{A_j\right\}_{j\geq 1}$  and corresponding points $\left\{x_{{A_j}}\right\}_{j\geq 1}$ that are required to be dense in  $supp(\Measure[]{B})$, the support of the measure $\Measure[]{B}$ (see line after Equation (2.4) on page 272 in~\cite{Diaconis:Ylvisaker:1979}). 
 It can be shown that a similar construction in the case of the $\ML$ distribution would lead to a $x_{_A}$ where $x_{_A}$ does not belong to $supp(\Measure[]{B})$, the Stiefel manifold. 
 Therefore,  the mentioned set of points $\left\{x_{{A_j}}\right\}_{j\geq 1}$ that are dense in $supp(\Measure[]{B})$ does not exist for the case of the $\ML$  distribution. 

Together, the two observations make it evident that  Theorem $1$ in~\citep{Diaconis:Ylvisaker:1979} is not applicable for constructing conjugate priors for the $\ML$ distribution. We would like to point out that the construction of the class of priors in~\cite{Hornik:2013} is based on a direct application of DY, which is not entirely applicable for the $\ML$-distribution. On the other hand,  the idea of constructing a conjugate prior on the natural parameter $F$ followed by a transformation, involves calculations of a complicated Jacobian term~\citep{Hornik:2013}. Hence the class of priors obtained via this transformation lacks interpretation of the corresponding hyperparameters.

\subsection{Two novel classes of Conjugate Priors}
\label{subsec:prior_novel}

Let $\Measure[]{M}$ denote the  normalized Haar measure on $\StiefelS$, $\Measure[]{V}$ denote the normalized Haar measure on $\SpaceV$, and  $\Measure[]{D}$ denote the Lebesgue measure on $\mathbb{R}_{+}^{p}$. For the parameters of the $\ML$-distribution, we define the prior density with respect to the product measure $\Measure[]{M} \times \Measure[]{D} \times \Measure[]{V}$ on the space $\StiefelS \times \Rplus^p \times \SpaceV$.
\begin{defn}
\label{defn:joint_prior}
The probability density function of the joint conjugate prior on the parameters $M, \bd$ and $V$ for the $\ML$ distribution is proportional to 
\begin{eqnarray}
g(M,\bd,V \,;\,\nu,\priorXzero) = \frac{etr\left( \nu\,VDM^T\priorXzero\right)}{ \left[_0 F_1 (\frac{n}{2}, \frac{D^2}{4})\right]^{\nu}},
\end{eqnarray}
as long as $ g(M,\bd,V \,;\,\nu,\priorXzero)$ is integrable.
Here $\nu > 0$ and $\priorXzero \in \Rnp$.  
\end{defn}
Henceforth, we refer to the joint distribution  corresponding  to the probability density function  in Definition~\ref{defn:joint_prior} as the joint conjugate prior distribution ($\JCPD$). We use the terminology,  joint conjugate prior class ($\JMDY$) when we use
\begin{eqnarray}\label{eq:JMDY:PriorClass}
(M, \bd, V) \sim \JCPD\left( \cdot \;;\nu,\priorXzero  \right),
\end{eqnarray}
as a prior distribution for the parameters of the $\ML$-distribution.
Although, the $\JMDY$ has some desirable properties (see Theorem~\ref{thm:DY_MDV_property} and Section~\ref{subsec:linearity_modal_param}), it may not be adequately flexible to incorporate prior knowledge about the parameters if the strength of prior belief is not uniform across the different parameters. For example, if a practitioner has strong prior belief for the values of $M$ but is not very certain about parameters $\bd$ and $V$, then $\JMDY$ may not be the optimal choice. Also, the class of joint prior defined in Definition~\ref{defn:joint_prior} corresponds to a dependent prior structure for the parameters $M$, $\bd$ and $V$. However, it is customary to use independent prior structure for parameters of curved exponential families~\citep{Casella:2002,Gelman:2014, Khare:2017}. Consequently, we also develop a class of conditional conjugate prior where we assume independent priors on the parameters $M$, $\bd$ and $V$. This class of priors are flexible enough to incorporate prior knowledge about the parameters even when the strength of prior belief differs across different parameters.

It is easy to see that the conditional conjugate priors for both $M$ and $V$ are $\ML$-distributions whereas the following definition is used to construct the conditional conjugate prior for $\bd$.

\begin{defn}
\label{defn:indep_prior}
The probability density function of the conditional conjugate prior for $\bd$ with respect to the Lebesgue measure on $\Rplus^p$ is proportional to 
\begin{eqnarray}
g(\bd \,;\,\nu,\BoEta, n) = \frac{ \exp(\nu\,\BoEta^T\bd)}{\left[_{0}F_1 \left(  \frac{n}{2}, \frac{D^2}{4}\right)\right]^{\nu}},
\end{eqnarray}
as long as $g(\bd \,;\,\nu,\BoEta,n)$ is integrable. Here $\nu > 0$, $\BoEta \in \mathbb{R}^p$ and $n\geq p$.
\end{defn}

Note that $g(\bd \,;\,\nu,\BoEta)$ is a function of $n$ as well.  However we do not vary $n$ anywhere in our construction, and thus we omit reference to $n$ in the notation for $g(\bd \,;\,\nu,\BoEta)$.

Henceforth we use the terminology, conditional conjugate prior distribution for  $\bd$ ($\CCPD$) to   refer to the probability distribution corresponding to the probability density function in Definition~\ref{defn:indep_prior}. We use the phrase  conditional conjugate prior class (\IMDY), to refer to the following structure of prior distributions 
\begin{eqnarray}\label{eq:IMDY:PriorClass}
M &\sim& \ML\left(\cdot;\, \hyparamM{M},\hyparamM{D},\hyparamM{V}  \right),\nonumber\\
  \bd &\sim& \CCPD\left(\cdot ;\, \nu,\BoEta\right),\nonumber\\
   V &\sim& \ML\left( \cdot;\,\hyparamV{M},\hyparamV{D},\hyparamV{V}\right),
\end{eqnarray}
where $M,\bd,V$ are assumed to be independent apriori.  As per Definitions ~\ref{defn:joint_prior} and \ref{defn:indep_prior}, the integrability of  the kernels mentioned in (3) and (5) are critical to   prove the propriety of the proposed class of priors. In light of this, Theorem~\ref{thm:DY_joint_prior} and Theorem~\ref{thm:DY_indep_prior}  provide conditions on $\nu, \priorXzero$ and $\BoEta$ for $g(M,\bd,V \,;\,\nu,\priorXzero)$ and $g(\bd \,;\,\nu,\BoEta)$ to be integrable, respectively. 
 \begin{theorem}
\label{thm:DY_joint_prior}
Let  $M\in \StiefelS$, $V \in \SpaceV$  and  $\bd\in \Rplus^p$. Let $\priorXzero \in \Rnp$ with $n\geq p$, then for any $\nu>0$,
\begin{enumerate}[(a)]
\item If $\normtwo{\priorXzero}<1$, then
\begin{eqnarray*}
\int_{\StiefelS}\int_{\SpaceV} \int_{\Rplus^p} g(M,\bd,V \,;\,\nu, \priorXzero)\;  \Measure[1]{D}\; \Measure[1]{V}\; \Measure[1]{M}<\infty,
\end{eqnarray*}
\item If $\normtwo{\priorXzero}>1$, then
\begin{eqnarray*}
\int_{\StiefelS}\int_{\SpaceV} \int_{\Rplus^p} g(M,\bd,V \,;\,\nu, \priorXzero)\; \Measure[1]{D}\; \Measure[1]{V}\; \Measure[1]{M} = \infty,
\end{eqnarray*}
\end{enumerate}
where $g(M,\bd,V ; \nu, \priorXzero)$ is defined in Definition~\ref{defn:joint_prior}.
\end{theorem}

The conditions mentioned in this theorem do not span all cases; we have not addressed the case where $\normtwo{\priorXzero} = 1$.
 As far as statistical inference for practical applications is concerned, we may not have to deal with the case where $\normtwo{\priorXzero}=1$ as the hyper-parameter selection procedure (see Section~\ref{sec:hyperparameter_singleML}) and posterior inference (even in the case of uniform improper prior, see Section~\ref{subsec:UniformImproperPrior} ) only involve cases with $\normtwo{\priorXzero}<1$. We therefore postpone further investigation into this case as a future research topic of theoretical interest. 

\begin{theorem}
\label{thm:DY_indep_prior}
Let  $\bd \in \mathbb{R}_{+}^p$,  $\BoEta=\left(\eta_1, \ldots, \eta_p\right)\in \mathbb{R}^p$ and $n$ be any integer with $n\geq p$. Then for any $\nu > 0$, 
\begin{eqnarray}
\int_{\mathbb{R}_{+}^{p}} g(\bd ; \nu, \BoEta, n)\; \Measure[1]{D}  < \infty, \nonumber
\end{eqnarray}
if and only if $\max\limits_{1 \leq j \leq p} \eta_j < 1 $, where $g(\bd ; \nu, \BoEta, n)$ is as defined in Definition~\ref{defn:indep_prior}.
\end{theorem}

We can  alternatively parametrize the $\CCPD$ class of densities by the following specification of the probability density function, 
$$f(\bd \,;\, \nu, \BoEta) \propto \frac{\exp\left(\sum_{j=1}^{p} \eta_j d_j\right)}{\left[ _0 F_1 (\frac{n}{2}, \frac{D^2}{4})\right]^{\nu}},$$
where $\max_{1 \leq j \leq p} \eta_{j} < \nu$.
In this parametrization, if we consider the parameter choices, $\nu = 0$ and $\boldsymbol{\beta}:=-\BoEta$, then the resulting  probability distribution corresponds to the {\it{Exponential}}  distribution with rate parameter  $\bbeta$.

{\paragraph{Properties of the $\CCPD$ and $\JCPD$ Class of Distributions}}


It is important to explore the properties for the $\CCPD$ and $\JCPD$ class of distributions in order to use them in an effective manner. Intuitive interpretations of the parameters $\nu, \BoEta, \priorXzero$ are desirable, for example, for hyper-parameter selection. Due to conjugacy, Bayesian analysis will lead to posterior distributions involving $\JCPD$ and $\CCPD$, and therefore, it is necessary to identify features that are required to develop practicable computation schemes for posterior inference.
The following four theorems establish some crucial properties of the $\CCPD$ and $\JCPD$ class of distributions.




\begin{theorem}
\label{thm:DY_D_property}
Let $\variableX \sim \mbox{\CCPD}(\cdot ; \nu, \BoEta)$ for $\nu>0$ and $\max_{1 \leq j \leq p} \eta_j<1$ where $\BoEta=\left(\eta_1,\ldots, \eta_p \right)$. Then  
\begin{enumerate}[(a)]
\item The distribution of $\bd$ is log-concave.
\item The distribution of $\bd$ has a unique mode if $\eta_j > 0$ for all $j=1, 2, \cdots, p$. The mode of the distribution is given by $\m = {\bf h}^{-1}(\BoEta)$, where the function $ {\bf h}(\bd)$  is defined as follows, ${\bf h}(\bd):= \left(h_1(\bd), h_2(\bd),\cdots ,h_p(\bd)\right)^T$ with $$h_j(\bd) := {\left(\frac{\partial }{\partial d_j}\,{}_0F_1 \left(\frac{n}{2}, \frac{D^2}{4} \right)\right)}/{_0F_1 \left(\frac{n}{2}, \frac{D^2}{4} \right)}.$$

\end{enumerate}
\end{theorem}

Notably, the mode of the distribution is characterized by the parameter $\BoEta$ and does not depend on the parameter $\nu$. The proof of the theorem  relies on a few nontrivial properties of $\hyp$, i.e., the hyper-geometric function of a matrix argument, that we have established in the supplementary material Section~\ref{s-Apndx:sec:hypProperties}. It is easy to see that the function ${\bf h}^{-1}$ is well defined as the function ${\bf h}$ is strictly increasing in all its coordinates. Even though subsequent theoretical developments are based on the formal definition and theoretical properties of ${\bf h}^{-1} $ and ${\bf h}$ functions, numerical computation of the functions are tricky. The evaluation of the functions depend on reliable computation of $\hyp$ and all its partial derivatives. In Section~\ref{sec:HYPComputation}, we provide a reliable and theoretically sound computation scheme for these functions.  

On a related note, it is well known that log-concave densities correspond to unimodal distributions if the sample space is the entire Euclidean space~\citep{Ibragimov:1956, Dharmadhikari:1988, Doss:2016}. However, the mode of the distribution may not necessarily be at a single point. Part(b) of Theorem~\ref{thm:DY_D_property} asserts that the $CCPD$ has a single point mode. Moreover, the sample space of $CCPD$ is $\bd\in\Rplus^p$, which merely encompasses the positive quadrant and not the whole of the $p$ dimensional Euclidean space. Hence general theories developed for $\mathbb{R}^p$ (or $\mathbb{R}$) do not apply. In fact, when $\eta_j \leq 0$, the density defined in Definition~\ref{defn:indep_prior} is decreasing as a function of $d_j$ on the set $\Rplus$ and the mode does not exist as $\Rplus$ does not contain the point $0$. In all, part(b) of Theorem~\ref{thm:DY_D_property} does not immediately follow from part(a) and requires additional effort to demonstrate.

In order to introduce the notion of ``concentration" for the $\CCPD$ class of distributions we require the concept of a level set. 
Let the unnormalized probability density function for the $\CCPD$ class of distributions, $g(\bx;\nu,\BoEta)$ (See Definition~\ref{defn:CCPD_conditional}), achieve its maximum value at $\m$ ( part(b) of Theorem~\ref{thm:DY_D_property} ensures that $\m$ is a unique point) and let 
\begin{eqnarray}
{\SetWithModePrime}_l = \left\{ \bx \in \Rplus^p: g(\bx;1,\BoEta)/g(\m;1,\BoEta) > l \right\}
\label{eq:levelSetDefinition}
\end{eqnarray}
be the level set of level $l$ containing the mode $\m$ where $0 \leq l < 1$. To define the level set we could have used $g(\bx;\nu_0,\BoEta)$ for any fixed value of $\nu_0 > 0$ instead of $g(\bx;1,\BoEta)$. However, without loss of generality, we choose $\nu_0 = 1$.

Let $P_{\nu }(\cdot ; \BoEta)$ denote the probability distribution function corresponding to the  $ \mbox{\CCPD}(\cdot ; \nu, \BoEta)$ distribution. According to Theorem\ref{thm:DY_D_property}, for a fixed $\BoEta\in \mathbb{R}^p$,  all distributions in the class $\{ P_{\nu }(\cdot ; \BoEta) : \nu>0\}$ have the mode located at the point $\m$. 

\begin{theorem}
\label{thm:DY_D_property_concentration}
Let $\variableX_{\nu} \sim \mbox{\CCPD}(\cdot ; \nu, \BoEta)$ for a fixed $\BoEta \in \mathbb{R}^p$  with $\m$ being the mode of the distribution. If   $P_{\nu}(\cdot ;\BoEta)$ denotes the probability distribution function corresponding to $\variableX_{\nu}$, then
\begin{enumerate}[(a)]
\item  $P_\nu(\SetWithMode_l; \BoEta)$ is an increasing function of $\nu$ for any level set $\SetWithMode_l$ with  $l \in (0,1)$,
\item  For any open set $\SetWithMode \subset \Rplus^p$ containing $\m$, $P_\nu(\bd \in \SetWithMode ;\BoEta)$ goes to $1$ as $\nu \to \infty$.
\end{enumerate}
\end{theorem}

The major impediment to proving Theorem~\ref{thm:DY_D_property_concentration} arises from the intractability of the normalizing constant of the $CCPD(\cdot ; \nu, \BoEta)$ distribution. Although involved, the proof essentially uses the log convexity of $\hyp$ to get around this intractability.

From Theorem~\ref{thm:DY_D_property_concentration}, it is clear that the parameter $\nu$ relates to the concentration of the probability around the mode of the distribution. Larger values of $\nu$ imply larger concentration of probability near the mode of the distribution.


\begin{defn}
\label{defn:modal_def}
In the context of the probability distribution  $\CCPD\left(\cdot \;;\; \BoEta, \nu \right)$, the parameters $\BoEta$ and  $\nu$ are labeled as the ``modal parameter" and the  ``concentration parameter", respectively. 
\end{defn}
%

\begin{figure}
\centering
\begin{tabular}{ccc}
\hspace{-.5 in}\includegraphics[scale=0.5]{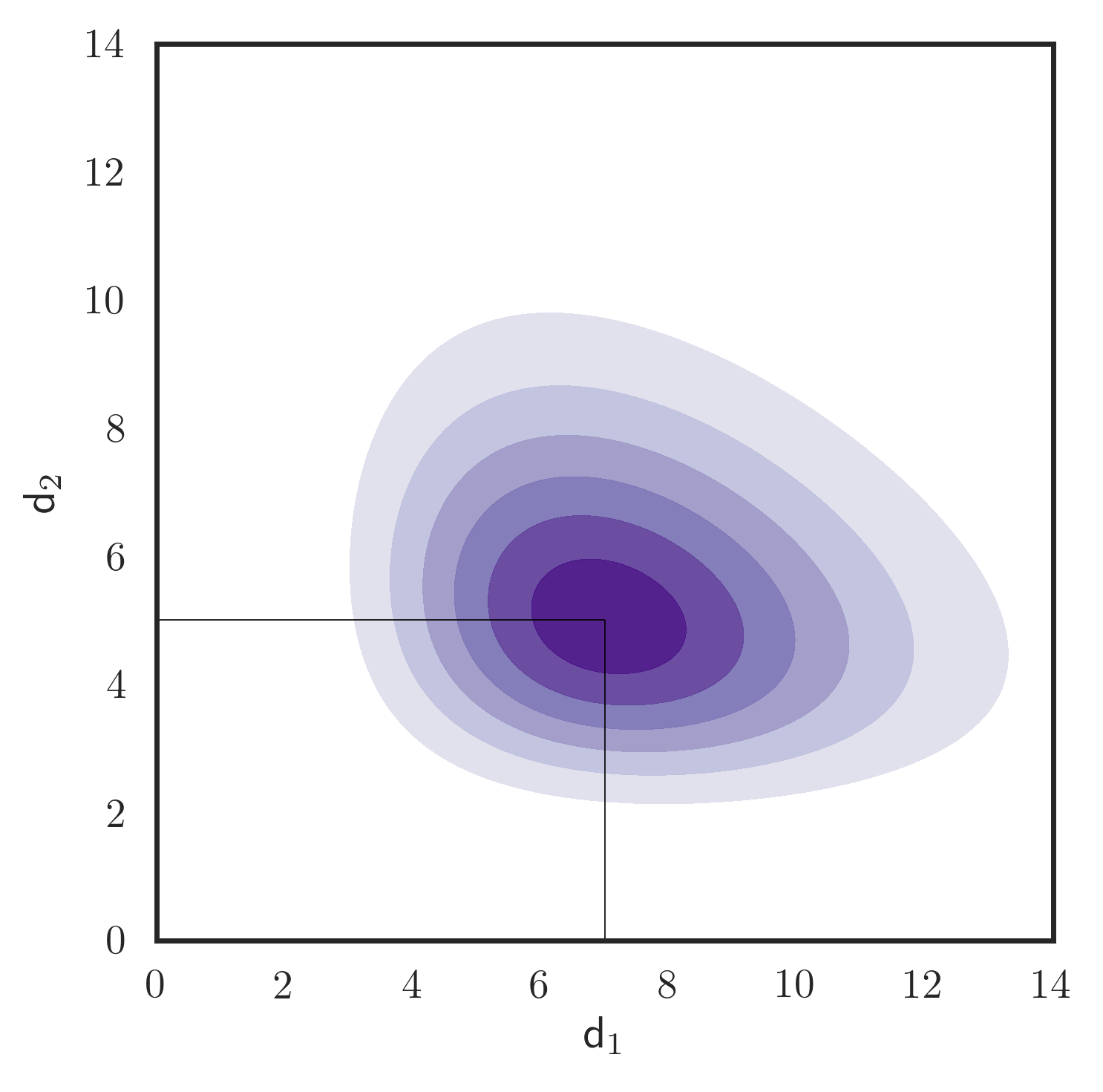}& 
\hspace{-.2  in}\includegraphics[scale=0.5]{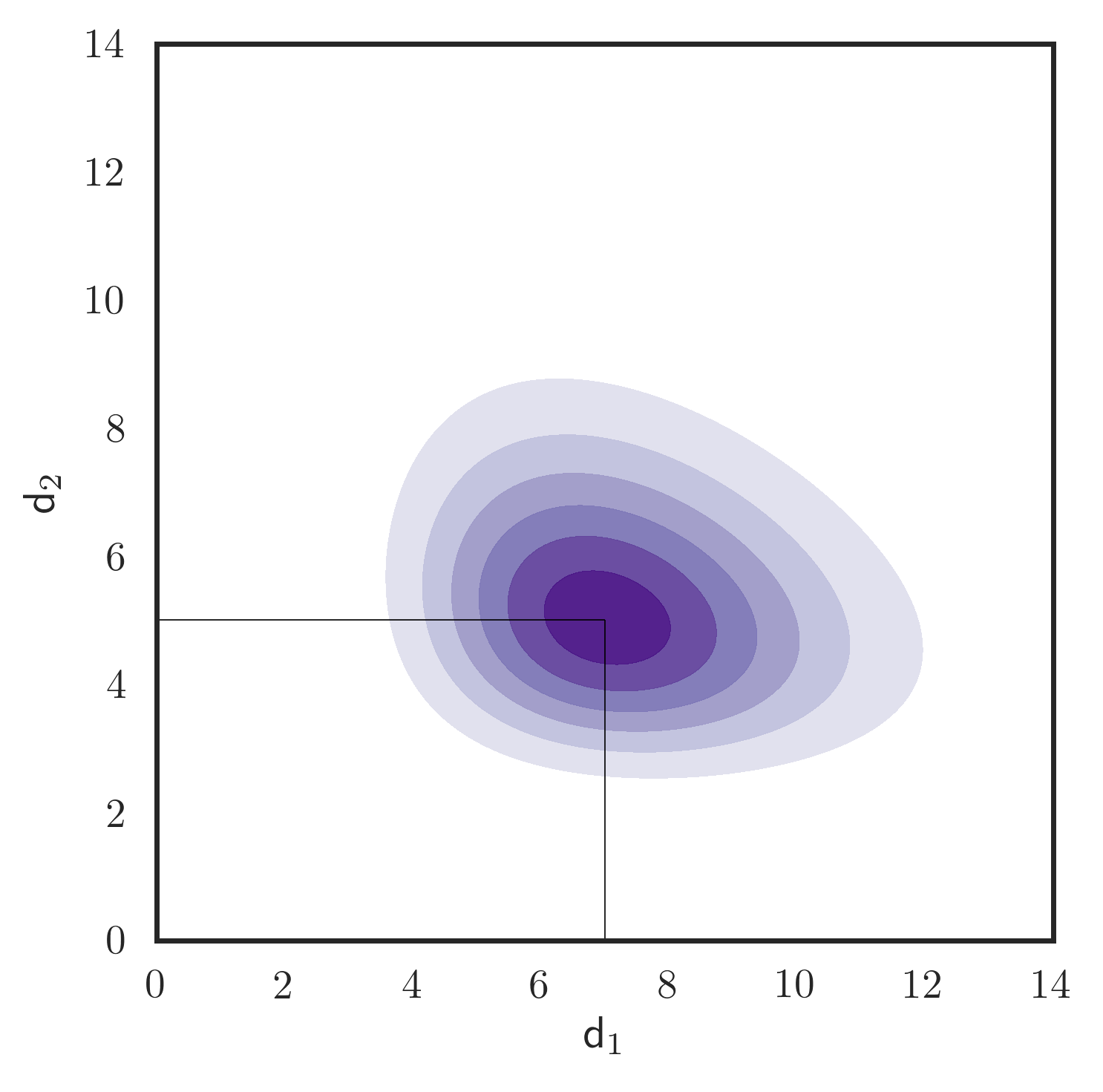} & 
\hspace{-.2 in}\includegraphics[scale=0.5]{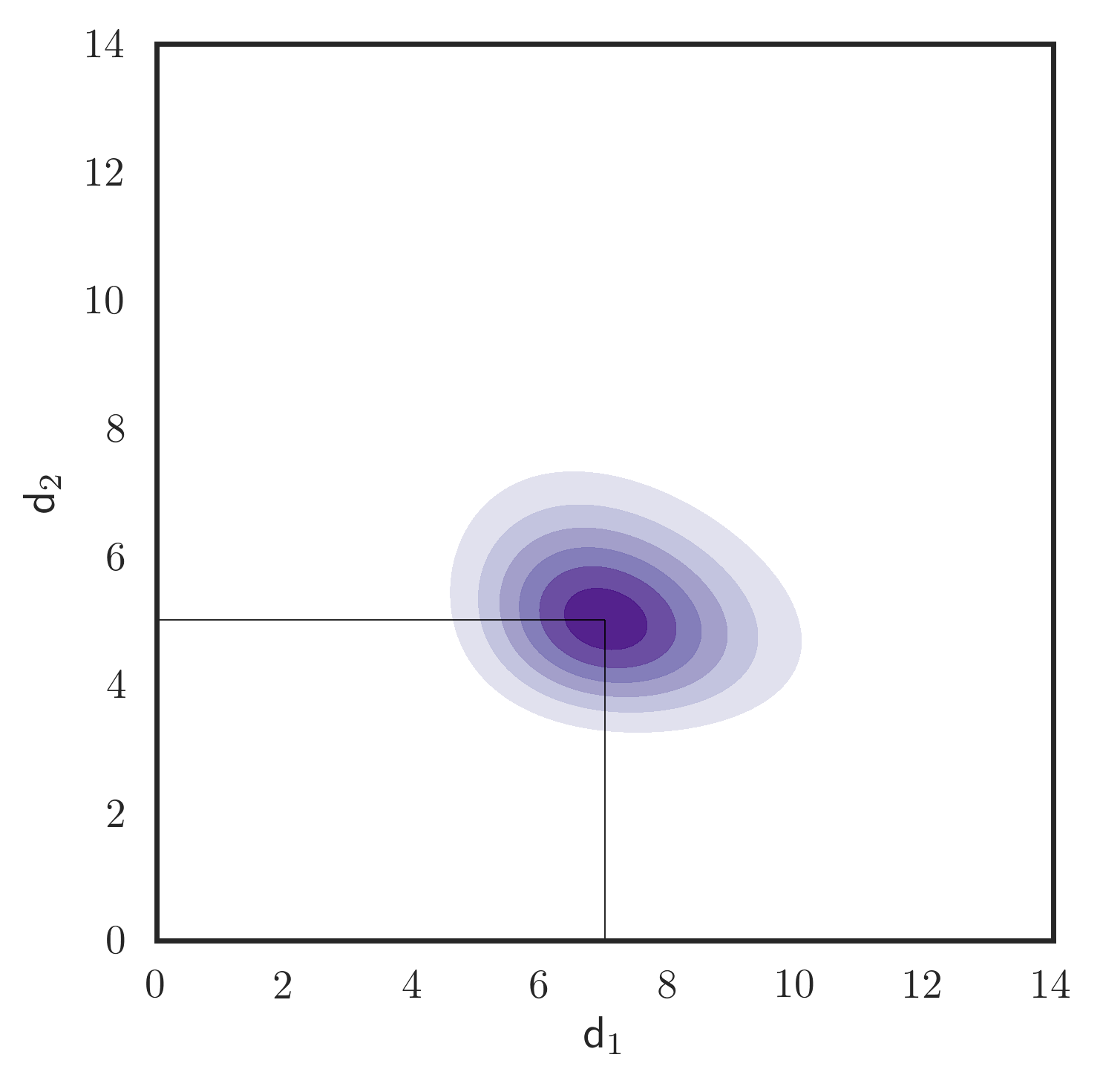}\\
\hspace{-.45 in}(a) $\nu = 10$ &
\hspace{-.1 in}(b) $\nu = 20$ &
(c) $\nu = 35$ \\
\end{tabular}
\caption{Density plots of $CCPD(\cdot \, ; \;\nu, \BoEta )$ for different values of $\nu$ where $\BoEta = (0.89,0.85)$. Mode of the distributions are located at the point $(7,5)$.}
\label{fig:prior_plots}
\end{figure}

In Figure~\ref{fig:prior_plots}, we display three contour plots of the $\CCPD(\cdot \;; \nu, \BoEta)$ distribution with ${\BoEta}=(0.85, 0.88)$. Note that the corresponding mode of the distribution is ${\bf h}^{-1}(0.85, 0.88)=(7, 5)$ for all three plots.  We can observe the implication of part~(b) of Theorem~\ref{thm:DY_D_property} as the ``center" of the distributions are the same. Contrastingly, it can be observed that the ``spread" of the distributions decrease as the value of the parameter $\nu$ increases, as implied by Theorem~\ref{thm:DY_D_property_concentration}.


\begin{theorem}
\label{thm:DY_MDV_property}
Let $(M,\variableX, V) \sim  \mbox{\JCPD}(\cdot ; \nu, \priorXzero)$ for some $\nu>0$ and $\normtwo{\priorXzero}<1$. If $\priorXzero = M_{\priorXzero}D_{\priorXzero} V_{\priorXzero}^T$ is the unique SVD of $\priorXzero$ with $\bd_{\priorXzero}$ being the diagonal elements of $D_{\priorXzero}$, then the unique mode of the distribution is given by 
$(M_{\priorXzero}, {\bf h}^{-1}(\bd_{\priorXzero}), V_{\priorXzero})$ 
where the function $\bd \rightarrow {\bf h}(\bd)$  is as defined in Theorem~\ref{thm:DY_D_property}.

\end{theorem}

Note that the mode of the distribution is characterized by the parameter $\priorXzero$ and does not depend on the parameter $\nu$. The proof of the theorem depends crucially on a strong result, a type of rearrangement inequality proved in~\cite{Kristof:1969}.

For the concentration characterization of $JCPD$, we define the level sets in the context of the $JCPD$ distribution. Let the unnormalized probability density function for the $\JCPD$ class of distributions, $g(M, \bd, V;\nu,\priorXzero)$, achieve its maximum value at the point $(\hat{M},\hat{\bd},\hat{V})$ ( see Theorem~\ref{thm:DY_MDV_property} ) and  
$$
{\LevelSetMDV}_l = \left\{ ( M, \bd, V) \in \StiefelS \times\Rplus^p \times \SpaceV: g(M,\bd,V;1,\priorXzero)/g( \hat{M},\hat{\bd},\hat{V}  ;1,\priorXzero) > l \right\}
$$ be the level set of level $l$ from some $l\in (0,1)$. The following theorem characterizes the concentration property of the $JCPD$ distribution. 

\begin{theorem}
Let  $(M,\bd, V) \sim \mbox{\JCPD}(\cdot ;\nu,\priorXzero)$, where $\normtwo{\priorXzero}<1$. If $P_{\nu}(\cdot\;;\;\priorXzero)$ denotes the probability distribution function corresponding to the distribution $ \mbox{\JCPD}(\cdot;\nu,\priorXzero)$, then
\begin{enumerate}[(a)]
\item  $P_\nu(\LevelSetMDV_l\;;\;\priorXzero)$ is a strictly increasing function of $\nu$ for any level set $\LevelSetMDV_l$ with  $l \in (0,1)$.
\item  For any open set $\SetWithModeMDV  \subset \StiefelS \times\Rplus^p \times \SpaceV$ containing the mode of the distribution, $P_\nu( \SetWithModeMDV \;;\;\priorXzero)$ tends to $1$ as $\nu \to \infty$.
\item The conditional distribution of $M$ given $(\bd,V)$ and $V$ given $(M,\bd)$ are $\ML$ distributions whereas the conditional distribution of $\bd$ given $(M,V)$ is a $\CCPD$ distribution. 
\end{enumerate}
\label{thm:DY_MDV_property:concentration}	
\end{theorem}
%

Parts (a) and (b) of the above theorem characterize the concentration whereas part(c) relates $CCPD$ to the $JCPD$ class of distributions. Part(c) also motivates the development of a sampling procedure for the $JCPD$ distribution. The proof of part(a) Theorem~\ref{thm:DY_MDV_property:concentration} is similar to that of the proof of Theorem~\ref{thm:DY_D_property_concentration}. The proof for part(b) of Theorem~\ref{thm:DY_MDV_property:concentration} is more involved and depends on several key results, including the rearrangement inequality by \citep{Kristof:1969}, the log convexity of $\hyp$, and the the fact that $g({\bf h}^{-1}(\BoEta)\; ; \; \nu, \BoEta))$, the value of the unnormalized $ CCPD$ density at the mode, is a strictly increasing function of the parameter $\BoEta$.

Note that unlike in the case of the $CCPD$ distribution, we do not attempt to establish the log concavity of $JCPD$, the reason being that the underlying probability space $\StiefelS \times\Rplus^p \times \SpaceV$ is non-convex. Nevertheless, it is evident that beyond a certain distance (based on a suitable metric on  $\StiefelS \times\Rplus^p \times \SpaceV$) the value of the density drops monotonically as one moves farther away from the center. 
Based on the characteristics of the parameters $\nu$ and $\priorXzero$ of the $\JCPD$ class of distributions,  we have the following definitions.

\begin{defn}
\label{defn:modal_def:JCPD}
The parameters $\priorXzero$ and $\nu$  in the distribution $\JCPD$ are labeled the ``modal" parameter and the ``concentration" parameter, respectively. 
\end{defn} 
 
Interestingly, both  distributions $CCPD$ and $JCPD$ are parameterized by two parameters, one controlling the center and the other characterizing the probability concentration around that center. One may therefore visualize the distributions in a fashion similar to that of the multivariate Normal distribution controlled by the mean and variance parameters. This intuitive understanding can help practitioners select hyper-parameter values when conducting a Bayesian analysis with the $CCPD$ and $JCPD$ distributions.
 
Thus far we have established properties of $CCPD$ and $JCPD$ that relate to basic features of these distributions. Additional properties, which are required for a MCMC sampling scheme, are developed in Section ~\ref{sec:posteriorDistribution}.
%

\section{Hyperparameter Selection Procedure}
\label{sec:hyperparameter_singleML}
\subsection{Informative Prior}
We now present procedures for the selection of hyperparameter values aimed at incorporating prior beliefs about the parameters $(M,\bd,V)$. Consider the scenario where a practitioner has the prior belief that the values for the parameters $M,\bd,V$ are close to $M_{belief},\bd_{belief},V_{belief}$, respectively. A standard approach to incorporating this prior knowledge is to select the hyper-parameter values in such a manner that the mode of the corresponding prior distribution becomes $M_{belief},\bd_{belief},V_{belief}$. In order to achieve this in the current context, we first  compute $\tilde{\BoEta}=h(\bd_{belief})$ where $h(\cdot)$ is defined in Equation~\ref{s-eq:priorModeEquation} in the supplementary material. Note that we always get a feasible $\tilde{\BoEta}$ for every real $\bd_{belief} \in \SpaceD$.

In the case of the $\IMDY$ class of priors, we choose $\BoEta=\tilde{\BoEta}$, $\hyparamM{M}=M_{belief}$, $\hyparamV{M}=V_{belief}$, $\hyparamM{V}=\iMat, \hyparamV{V}=\iMat$  in the Equation~\ref{eq:IMDY:PriorClass}. Theorem~\ref{thm:DY_D_property} guarantees that the above hyper-parameter specifications yields a prior distribution that has mode at  $(M_{belief}, \bd_{belief}, V_{belief})$. From Theorem~\ref{thm:DY_D_property}, we also see that larger values of the hyper-parameter $\nu$ lead to larger concentration of the prior probability around the mode. The hyper-parameters $\hyparamM{D}$ and $\hyparamV{D}$ play a similar role for the $\ML$ distribution.  Hence the hyper parameters $\nu, \hyparamM{D}$ and $\hyparamV{D}$  are chosen to have larger values in case the practitioner has a higher confidence in the prior belief.

In the case of the $\JMDY$ class of priors, we apply Theorem~\ref{thm:DY_MDV_property} to construct $\JCPD$ (see Equation~\ref{eq:JMDY:PriorClass}) with mode at $M_{belief},\bd_{belief},V_{belief}$. In particular, we set ${\Psi}=M_{belief}D_{\tilde{\BoEta}}(V_{belief})^T$ where $D_{\tilde{\BoEta}}$ is the diagonal matrix with diagonal elements ${\tilde{\BoEta}}=h(\bd_{belief})$. Using the concentration characterization described in Theorem~\ref{thm:DY_MDV_property}, the practitioner may choose the value of the hyper-parameter $\nu$ appropriately, where a larger value for the parameter $\nu$ implies greater confidence in the prior belief. 

It is noteworthy that for both the $\JMDY$ and $\IMDY$ class of priors,  there is an intimate connection between the sample size and the interpretation of the hyper-parameter $\nu$. As a heuristic one may envisage $\nu$ as incorporating ``information" equivalent to $\nu$ many historic observations of the model.

\subsection{Uniform improper prior} 
In the case where the practitioner does not have a prior belief about the parameter values, an automatic procedure for hyper-parameter selection can be helpful. In this and the next subsection, we  discuss two automatic procedures to select the values of the hyper-parameters.
 In the absence of prior information, usage of uniform prior is common in the literature. In the context of the current model, for the $\JMDY$ and $\IMDY$ class of distributions, the prior for the parameters $(M,\bd, V)$, is called a  uniform  prior if
$$ g(M,\bd, V; \nu,\Psi) \propto 1  \text{ and }$$
$$ \MLDensity(M\;;  \hyparamM{M},\hyparamM{D}, \hyparamM{V})g(\bd; \nu,\BoEta) \MLDensity(V\;;\;\hyparamV{M},\hyparamV{D}, \hyparamV{V}) \propto 1 .$$ 
 Both  classes of priors $\JMDY$ and $\IMDY$ are flexible enough to accommodate a uniform  prior. For $\JMDY$, this can be achieved by setting $\nu = 0$ in Equation~\ref{eq:JMDY:PriorClass}. Correspondingly, for the $\IMDY$ class, the uniform  prior can be constructed by choosing $\nu=0$, $\hyparamM{D}={\bf 0}$ and $\hyparamV{D}={\bf 0}$ in Equation~\ref{eq:IMDY:PriorClass}. Note that the resulting uniform prior is improper in nature as the above choices of hyper parameters do not lead to a proper probability distribution. Hence, it is necessary to check the propriety of the resulting posterior (see Section~\ref{subsec:UniformImproperPrior} for more details).

\subsection{Empirical prior} 
Another widely used automatic method is to use empirical information contained in the data to select appropriate values of the hyper-parameters. Let $W_1, W_2, \ldots W_N$ be independent and identically distributed samples drawn from $\ML( \cdot \;; M, \bd, V)$. Consider the sample mean, $\overline{W}=(\sum_{i=1}^{N}W_i)/N$. Let the unique SVD of the sample mean be $\overline{W}= M_{\overline{W}} D_{\overline{W}} V_{\overline{W}}$. Construct candidate values $M_{belief}=M_{\overline{W}}$, $V_{belief}=V_{\overline{W}}$ and $\tilde{\BoEta}$ as the diagonal elements of  $ D_{\overline{W}}$. One can set $\Psi=\overline{W}$ as the hyper-parameter in the case of the $\JMDY$ prior.  In the case of the $\IMDY$ class of priors, one can choose $\BoEta=\tilde{\BoEta}$, and for the hyper-parameters related to $M$ and $V$,  apply the same procedure as discussed previously in this section. For both classes of priors, a value for $\nu$ that is less than or equal to $10$ percent of the  sample size $N$, is recommended.

\begin{example}
Let the practitioner have the following prior belief for the values of the parameters $M, \bd, V$,
$$M_{belief}=
\begin{bmatrix}
    1       & 0 \\
    0       & 1  \\
    0       & 0 
\end{bmatrix},\;
\bd_{belief}=
\begin{bmatrix}
7\\
5
\end{bmatrix},\;
V_{belief}=
\begin{bmatrix}
    1       & 0 \\
    0       & 1 
\end{bmatrix} 
.$$
As described previously in this section, we can compute $\tilde{\BoEta}={\bf h}( 7,5)=(0.89, 0.85)$. Hence, for the $\JMDY$ class of priors,  we choose the hyper-parameter values 
$$\tilde{\priorXzero}= \begin{bmatrix}
    1       & 0 \\
    0       & 1  \\
    0       & 0 
\end{bmatrix} 
\begin{bmatrix}
    0.89       & 0 \\
    0       & 0.85
\end{bmatrix} 
 \begin{bmatrix}
    1       & 0 \\
    0       & 1 
\end{bmatrix}^T  = \begin{bmatrix}
    0.89       & 0 \\
    0       & 0.85  \\
    0       & 0 
\end{bmatrix}, $$ 
 to ensure that $\JCPD(\cdot \; ; \tilde{\priorXzero} , \nu)$ has mode at $M_{belief}, \bd_{belief}, V_{belief}$ for all values of $\nu>0$.  The value of the hyper-parameter $\nu$ should be chosen according to the strength of the prior belief.  In  Figure~\ref{fig:prior_plots}, we display the resulting conditional distribution for $\bd$ given $M,V$. Figure~\ref{fig:prior_plots} shows that the ``center" of the distribution is located at $(7, 5)$. Figure~\ref{fig:prior_plots} also displays the ``spread" of the distribution around the mode when using $\nu = 10$, $\nu = 20$ and $\nu=35$.
\end{example}

\section{Properties of Posterior}
\label{sec:posteriorProperties}

The derivation of the posterior distributions for the $\JMDY$ and $\IMDY$ class of priors is straightforward since they were built with conjugacy in mind, which then entails that the posterior distributions lie in the corresponding classes. However, inference for the resulting posterior distributions is challenging because not only are the normalizing constants intractable for both the $JCPD$ and $CCPD$ distributions, but also, the unnormalized version of the corresponding density functions involve $\hyp$. We first focus our attention on developing properties of the posterior distribution when involving $\JMDY$ and $\IMDY$ priors. In particular, we  derive explicit forms of the posterior conditionals under different prior settings, the linearity of the posterior mode parameters and the strong consistency of the posterior mode.

\subsection{Posterior conditionals}
\label{sec:posteriorDistribution}

Let $W_1, W_2, \ldots W_N$ be independent and identically distributed samples drawn from $\ML( \cdot \;; M, \bd, V)$. Let $\overline{W} = \sum_{i=1}^N W_i/N$. The likelihood of the data is 
\begin{eqnarray}\label{eq:likelihood}
 \prod_{i=1}^N \frac{etr(VDM^T W_i)}{{}_0F_1( \frac{n}{2}, \frac{D^2}{4})}.
\end{eqnarray}
First, let us assume a $\JCPD$  prior with parameters $\nu$ and $\priorXzero$. 
 Theorem~\ref{thm:DY_MDV_property} not only implies that the posterior has a unique mode, but also provides an expression for the mode.
Furthermore, we see that the corresponding  posterior distribution is $\JCPD $ with concentration $({\nu+N})$ and posterior modal parameter $\widehat{\priorXzero}_{N}= \left(\frac{\nu}{\nu+N}\priorXzero + \frac{N}{\nu+N} \overline{W}\right).$  Let $\hat{\BoEta}_{\priorXzero_N}$ be the diagonal elements of the diagonal matrix $\hat{D}_{\priorXzero_N}$, where $\widehat{\priorXzero}_{N} = \hat{M}_N \hat{D}_{\priorXzero_N}\hat{V}_N$ is the unique SVD for $\widehat{\priorXzero}_{N}$.
 From Theorem~\ref{thm:DY_MDV_property:concentration},  it follows that the full posterior conditionals for the parameters $M, \bd, V$ are  $\ML$, $\CCPD$ and $\ML$ distributions, respectively.
 
In Section \ref{subsubsec:post_comp} we shall use these results to construct a  Gibbs algorithm. A part of the Gibbs scheme would require sampling from the relevant $CCPD$ distribution, which we propose to implement by simulating from the full conditional distribution of each of the components of $\bd$ given the rest, when $\bd \sim \CCPD(\cdot \; ;  \nu,\BoEta)$. To refer to this conditional distribution in subsequent text, we have the following definition. 

\begin{defn}\label{defn:CCPD_conditional}
Let  $\nu>0$, $\boldsymbol{\varpi} \in \Rplus^{p-1}$ and $\BoEta\in \Rplus^{p}$ with $\max_{1\leq j\leq p}\eta_j<1$. A random variable is defined to be distributed as $\CCPD^{\star}_j \left(\cdot\; ;\boldsymbol{\varpi}, \nu ,\BoEta\right)$, if the corresponding probability density function (with respect to the Lebesgue measure on $\mathbb{R}$) is proportional to 
\begin{eqnarray}
g_j( x ;\,\boldsymbol{\varpi},\nu,\BoEta) =\frac{ \exp(\nu\,\eta_j x)}{\left[_{0}F_1 \left(  \frac{n}{2}, \frac{(\Delta(x))^2}{4}\right)\right]^{\nu}},\nonumber
\end{eqnarray}  
where $\Delta(x)$ is a diagonal matrix with diagonal elements $(x,\boldsymbol{\varpi})\in \Rplus^{p}.$
\end{defn}

Let $\bd=(d_1, \ldots, d_p)$ be a random vector with $\bd \sim \CCPD\left(\cdot \; ; \nu,\BoEta \right)$ for some $\max_{1\leq j\leq p}\eta_j<1, \nu>0$. Let $\bd^{(-j)}$ be the vector containing all but the $j$-th component of the vector $\bd$. Then the conditional distribution  of $d_j$ given $ \bd^{(-j)}$ is   $\CCPD_{j}^{\star}(\cdot \; ;\bd^{(-j)}, \nu, \BoEta)$, i.e., 
$$ d_j \mid \bd^{(-j)} \sim \CCPD_{j}^{\star}(\cdot \; ;\bd^{(-j)}, \nu,\BoEta).$$

Now, since the conditional posterior of $\bd$ was shown to be $\CCPD$, the conditional posterior distribution of $d_j \mid \bd^{(-j)}, M, V, \AllData$ follows a $\CCPD_{j}^{\star}$ distribution. 
 

\paragraph{}
In the case of a Bayesian analysis with a $\IMDY$ prior, Equation~\ref{eq:IMDY:PriorClass} and~\ref{eq:likelihood} determine the corresponding  posterior distribution to be proportional to 

\begin{eqnarray}
 \frac{etr\left(\left(V\,D\,M^T\right)\,N\,\overline{W} + G^0\,M + H^0\,V\right)} {{}_0F_1(\frac{n}{2};D^2/4)^{\nu + N}}  
\exp(\nu\,\BoEta^T{\bd}),
\label{eq:post_mix}
\end{eqnarray}
 where $G^0 = \hyparamM{V}\,\hyparamM{D}\,{(\hyparamM{M})}^T${and} $ H^0 = \hyparamV{V}\,\hyparamV{D}\,{(\hyparamV{M})}^T$.
The conditional probability density for the posterior distribution of $\bd$ given $M$, $V$, ${\{W_i\}}_{i=1}^N$ is proportional to 
\begin{eqnarray}
 \frac{ \exp\left((\nu+N)\;{\left(\frac{\nu}{\nu+N}\BoEta + \frac{N}{\nu+N}\BoEta_{\overline{W}}\right)}^T\bd \right)}{\left[_{0} F_1 \left(  \frac{n}{2}, \frac{D^2}{4}\right)\right]^{\nu+N}}, \;\;
\label{eq:posterior_cond_IMDY}
\end{eqnarray}
where $\BoEta_{\overline{W}} = \left(Y_{1,1},\cdots, Y_{p,p}\right)$ with $Y = M^T\overline{W} V $. It follows that the conditional posterior distribution of $\bd$  given $M, V , \AllData$ is $\CCPD(\cdot \; ; \; \hat{\nu}_N, \hat{\eta}_N)$ where $\hat{\nu}_N=\nu+N$ and $ \hat{\eta}_N=\left(\frac{\nu}{\nu+N}\BoEta + \frac{N}{\nu+N}\BoEta_{\overline{W}}\right) $. The conditional posterior distributions $M\mid \bd, V,{\{W_i\}}_{i=1}^N $ and  $V\mid \bd, M,{\{W_i\}}_{i=1}^N $ are $\ML$ distributions. 


\subsection{Linearity of posterior modal parameter} 
\label{subsec:linearity_modal_param}
We observe that the posterior  modal parameter is a convex combination of the prior modal parameter and the sample mean when applying the $\JMDY$ class of priors. In particular, from Section~\ref{sec:posteriorDistribution} we get  $$\hat{\Psi}_N=\left(\frac{\nu}{\nu+N}\priorXzero + \frac{N}{\nu+N} \overline{W}\right).$$  

In a similar fashion, we observe from Equation~\ref{eq:posterior_cond_IMDY} that the modal parameter for the conditional posterior distribution of $\bd$ given $M,V, \AllData$ is a convex combination of the prior modal parameter and an appropriate statistic of the sample mean. We should point out here that the posterior linearity of the {\it natural parameter} of an exponential family distribution directly follows from \cite{Diaconis:Ylvisaker:1979}. However, in our parametrization, the ML density is a curved exponential family of its parameters  and posterior linearity appears to hold for the  ``modal parameter".

\subsection{Posterior propriety when using uniform improper prior}\label{subsec:UniformImproperPrior}

In the case where a uniform improper prior is used, the corresponding posterior is proportional to
\begin{eqnarray}\label{eq:Improper:Posterior}
\frac{etr\left(N\;\,VDM^T\overline{W}\right)}{ \left[_0 F_1 (\frac{n}{2}, \frac{D^2}{4})\right]^{N}},
\end{eqnarray}
where $\overline{W}=\frac{1}{N}\sum_{i=1}^NW_i$ (see Equation~\ref{eq:likelihood}). It follows from Theorem~\ref{thm:DY_joint_prior} that the function in Equation~\ref{eq:Improper:Posterior} leads to a proper distribution, $\JCPD(\cdot\; ;\;N, \overline{W})$,  if $\normtwo{ \overline{W}}<1$. The following theorem outlines the conditions under which  $ \normtwo{ \overline{W}}<1$.
 
 \begin{theorem}\label{thm:Jupp:mardia:1979}
 Let $W_1, \ldots, W_N$ be independent and identically distributed samples from an $\ML$-distribution on the space $\StiefelS$. If  
\begin{enumerate}[(a)] 
 \item $N\geq2$, $p<n$
 \item $N\geq 3$ , $p=n\geq 3$, 
 \end{enumerate}
then $\normtwo{\overline{W}}<1$ with probability 1,  where $\overline{W}=\frac{1}{N}\sum_{i=1}^NW_i$.
 \end{theorem}

\subsection{Strong consistency of the posterior mode}
In the case where we use a $ \JCPD(\cdot;\nu, \Psi)$ prior for Bayesian analysis of the data $\AllData$, the corresponding posterior distribution is a $\JCPD $ with concentration ${\nu+N}$ and posterior modal parameter $\widehat{\priorXzero}_{N}= \left(\frac{\nu}{\nu+N}\priorXzero + \frac{N}{\nu+N} \overline{W}\right)$ (See Section~\ref{sec:posteriorDistribution}). Let $\widehat{\priorXzero}_{N} = M_{\priorXzero}D_{\priorXzero} V_{\priorXzero}^T$ be the unique SVD of $\widehat{\priorXzero}_{N}$ with $\bd_{\priorXzero}$ being the diagonal elements of $D_{\priorXzero}$. Then from Theorem~\ref{thm:DY_MDV_property}, the unique mode of the distribution is given by $(\hat{M}_N,\hat{\bd}_N,\hat{V}_N)$ where 
$$\hat{M}_N=M_{\priorXzero}, \hat{\bd}_N= {\bf h}^{-1}(\bd_{\priorXzero}) \text{  and } \hat{V}_N=V_{\priorXzero}.$$ 
The form of the function $ {\bf h}(\bd)$  is provided in Theorem~\ref{thm:DY_D_property}. The nontrivial aspect of finding the posterior mode is the computation of the function ${\bf h}^{-1}(\bd_{\priorXzero})$. In our applications, we use a Newton-Raphson procedure to obtain ${\bf h}^{-1}(\bd_{\priorXzero})$ numerically. We use large and small argument approximations for $\hyp$ ( See \cite{Jupp:1979}) to initialize the Newton-Raphson algorithm for faster convergence. Note that the success of the Newton-Raphson procedure here depends on the efficient computation of $\hyp$ and its partial derivatives. In Section~\ref{sec:HYPComputation}, we provide a method to compute these functions reliably.

The following theorem demonstrates that the mode of the posterior distribution is a strongly consistent estimator for the parameters $M,\bd,V$.
 
\begin{theorem}\label{thm:Consistency:mode}
Let $W_1, \ldots, W_N$ be independent and identically distributed samples from $\ML(\cdot\; ; \; M,\bd,V)$. Let $\hat{M}_N,\hat{\bd}_N$ and $\hat{V}_N$ be the posterior mode when a $\JMDY$ prior is used. The statistic $\hat{M}_N,\hat{D}_N$ and $\hat{V}_N$ are consistent estimators for the parameters $M,D$ and $V$. Moreover
$$ (\hat{M}_N,\hat{\bd}_N, \hat{V}_N ) \stackrel{a.s.}{\longrightarrow} (M,\bd,V)  \text{ as } N\longrightarrow \infty,$$
where a.s. stands for almost sure convergence.
\end{theorem}

\section{MCMC sampling from the Posterior}\label{subsubsec:post_comp}

Apart from finding the posterior mode, a wide range of statistical inference procedures including point estimation, interval estimation (see Section~\ref{sec:real_data_app}) and statistical decision making (see Section~\ref{sec:real_data_app}) can be performed with the help of samples from the posterior distribution. For the $JCPD$ and $CCPD$ classes of distributions, neither is it possible to find the posterior mean estimate via integration, nor can we directly generate i.i.d. samples from the distributions. We therefore develop procedures to generate MCMC samples \yy using a Gibbs sampling procedure, which  requires the results on posterior conditionals stated in Section \ref{sec:posteriorDistribution}. \jj
 
\yy It follows from \jj Theorem ~\ref{thm:DY_MDV_property:concentration} and Section 5.1  that  \yy under $JCPD$ prior \jj the conditional distribution of $M$ given $\bd, V$ and the conditional distribution  of $V$ given $M, \bd$ are $\ML$ distributions, while the conditional distribution of $\bd $ given $M,V$ is $CCPD$.  Consequently, the conditional distribution of $d_j \mid \bd^{(-j)}, M, V, \AllData$ follows a $\CCPD_{j}^{\star}$ distribution (see Definition~\ref{defn:CCPD_conditional}).
Also, let us assume that the unique SVD for $\hat{\nu}_N\,{(\hat{\priorXzero}_N V D)} = M_{\widehat{\priorXzero}}^M D_{\widehat{\priorXzero}} ^M {(V_{\widehat{\priorXzero}}^M)}^T$ and for
$\hat{\nu}_N\,{(\hat{\priorXzero}_N^T M D)} =  M_{\widehat{\priorXzero}}^V D_{\widehat{\priorXzero}} ^V {(V_{\widehat{\priorXzero}}^V)}^T$. Also, let us denote the  vector containing the diagonal element of the matrix $M^T \hat{\priorXzero}_N V$ to be $\BoEta_{\widehat{\priorXzero}}$. \yy Based on the above discussion, we can now \jj describe  the algorithm as follows.
\begin{algorithm}[H]
    \begin{algorithmic}[1]
    \State{ Sample $M \mid \bd,V,\AllData \sim \ML\left( \cdot \; ; \; M_{\widehat{\priorXzero}}^M ,  {\bd}_{\widehat{\priorXzero}} ^M , {V_{\widehat{\priorXzero}}^M} \right)$,}
            \State { Sample $d_j  \mid \bd^{(-j)}  M,V,\AllData \sim \CCPD_{j}^{\star}\left( \cdot \; ; \bd^{(-j)}, \hat{\nu}_N ,\BoEta_{\widehat{\priorXzero}}  \;  \right)$ for $j =1\ldots p,$}
        \State { Sample $V \mid \bd,V,\AllData \sim \ML\left( \cdot \; ; \; M_{\widehat{\priorXzero}}^V,  {\bd}_{\widehat{\priorXzero}}^V, {V_{\widehat{\priorXzero}}^V} \right)$.}
    \end{algorithmic}
    \caption{Gibbs sampling algorithm to sample from posterior when using $\JMDY$ prior}
\end{algorithm}

If instead we use a $\IMDY$ prior, (see Equation~\ref{eq:IMDY:PriorClass}) for Bayesian analysis of the data, then the  full conditional distribution of $M,\bd,V$ are $\ML$, $\CCPD$ and $\ML$ distributions, respectively. The steps involved in the Gibbs sampling Markov chain are then as follows.

\begin{algorithm}[H]
    \begin{algorithmic}[1]
    \State{ Sample $M \mid \bd,V,\AllData \sim \ML\left(\cdot\,;\,S^M_G, S^D_G, S^V_G\right),$}
            \State { Sample $d_j \mid \bd^{(-j)},M,V,\AllData \sim \CCPD_{j}^{\star}\left( \cdot \; ; \bd^{(-j)},\hat{\nu}_N ,\hat{\eta}_N,\;  \right)$ for $j =1, \ldots p,$}
        \State { Sample $V \mid M,\bd,\AllData \sim \ML\left(\cdot\,;\,S^M_H, S^D_H, S^V_H\right),$}
    \end{algorithmic}
    \caption{Gibbs sampling algorithm to sample from posterior when using $\IMDY$ prior}
\end{algorithm}
where $\hat{\nu}_N$, $\hat{\eta}_N$ are defined in Equation~\ref{eq:posterior_cond_IMDY}  and ($S^M_G, S^D_G, S^V_G$) ,  ($S^M_H, S^D_H, S^V_H$) are the unique SVD of the matrices $(DV^T\,N\overline{W}^T + G^0)$ and $(DV^T\,N\overline{W}^T + H^0)$, respectively.\\

\yy To implement the above algorithms we  need to sample from the $\ML$ and $CCPD$ distributions. For the former, \jj we use the procedure developed in ~\citep{Hoff:2009} to sample from the $\ML$ distributions. \yy Sampling from $\CCPD_{j}^{\star}$ is much more involved and is explained in detail in the next subsection. The following result provides some theoretical guarantees that shall be useful for this specific sampler.\jj 

\begin{theorem}
\label{lem:right_tail_prob_bound}
Let $\bd \sim \CCPD(\cdot ; \nu, \BoEta)$ for some $\nu>0$ and $\BoEta = \left(\eta_1,\ldots, \eta_p \right)$ where $\max_{1 \leq j \leq p} \eta_j < 1$. Let $g_1(\cdot\,;\,\bd^{(-1)},\nu,\BoEta)$ denote the unnormalized density corresponding to $\CCPD_{1}^{\star}(\cdot \; ;\bd^{(-1)}, \nu,\BoEta)$, the conditional distribution of  $d_1$ given $(d_2, \ldots , d_p)$. 
\begin{enumerate}[(a)]
\item The probability density function corresponding to  $\CCPD_{1}^{\star}(\cdot \; ;\bd^{(-1)}, \nu,\BoEta)$ is log-concave on the support $\Rplus$.
\item If $0<\eta_1<1$, the distribution  $\CCPD_{1}^{\star}(\cdot \; ;\bd^{(-1)}, \nu,\BoEta)$ is unimodal and the mode of the distribution is given by $m$ where $h_1(m)=\eta_1$. If $\eta_1\leq 0$ then the probability density is strictly decreasing on $\Rplus$.  
\item If  $B>m$ is such that  $\frac{g_1(B;\,\bd^{(-1)},\nu,\BoEta)}{g_1(m;\,\bd^{(-1)},\nu,\BoEta)}<\epsilon$ for some $\epsilon>0$, then $P(d_1 > B \mid d_2, \ldots, d_p)< \epsilon$,
\item Let $M_{crit}$ be any positive number, then for all $d_1>M_{crit}$,
\begin{eqnarray}
g_1(d_1 \,;\,\bd^{(-1)},\nu,\BoEta)&\leq & K^{\dagger}_{n,p,M_{crit}} \;d_1^{\nu(n-1)/2}\; exp(\;-\nu(1-\eta_1)\,d_1),\nonumber\\
\end{eqnarray}
where $$ K^{\dagger}_{n,p,M_{crit}}= \left[ \frac{({p /4})^{\frac{n/2-1}{2}})}{\Gamma(n/2)\left\{  \sqrt{M_{cric}}\;\;e^{- M_{crit}}\;I_{n/2-1}(M_{crit})\right\} } \right]^{\nu}.$$
\end{enumerate}
\end{theorem}

Even though parts (a) and (b) of the above theorem follow immediately from Theorem~\ref{thm:DY_D_property}, they are included here for completeness; all the properties play a crucial role in the construction of the sampling technique for $CCPD^{\star}_j$. The proof of part(c) is essentially an implication of the fact that the right tail of the distribution decays at an exponential rate. To show part(d) we have developed a nontrivial lower bound for  $\hyp$. 

\begin{rmrk}\label{remark:KDaggar}
The constant  $K^{\dagger}_{n,p,M_{crit}}$ in  part(d) of Theorem~\ref{lem:right_tail_prob_bound} converges to a finite constant as  $M_{crit}$ approaches infinity. It follows from the properties of the Bessel function that  $$\lim_{M_{crit}\rightarrow \infty} \sqrt{M_{crit}} e^{-M_{crit}}I_{a-1}(M_{crit})=\frac{1}{\sqrt{2\pi}}$$ for all $a\geq \frac{3}{2}.$ Hence for larger values of $M_{crit}$, the value of $K^{\dagger}_{n,p,M_{crit}}$ approaches $\left[ \frac{\sqrt{2\pi}({p /4})^{\frac{n/2-1}{2}})}{\Gamma(n/2) } \right]^{\nu}$, a nonzero finite constant depending on $n,p,\nu$.
\end{rmrk}

Note that the ratio  ${g_1(B;\,\bd^{(-1)},\nu,\BoEta)}/{g_1(m;\,\bd^{(-1)},\nu,\BoEta)}$, mentioned in part(c), is free of the intractable normalizing constants of the distribution. Therefore, the numerical computation of the ratio is possible as long as we can compute the corresponding $\hyp$.
Using Theorem~\ref{lem:right_tail_prob_bound} , we develop  an  accept-reject sampling algorithm that can generate samples from $CCPD_j^{\star}$ with high acceptance probability. The detailed construction of the sampler is provided next. We conclude this section with a description of an efficient procedure for computing the $\hyp$ constant.


\subsection{\bf Efficient rejection sampler for the $\CCPD^{\star}_j$ distribution}\label{sec:rejectionSampler} 
\begin{figure}[!htb]
\centering
\begin{tabular}{cc}
\includegraphics[scale=0.5]{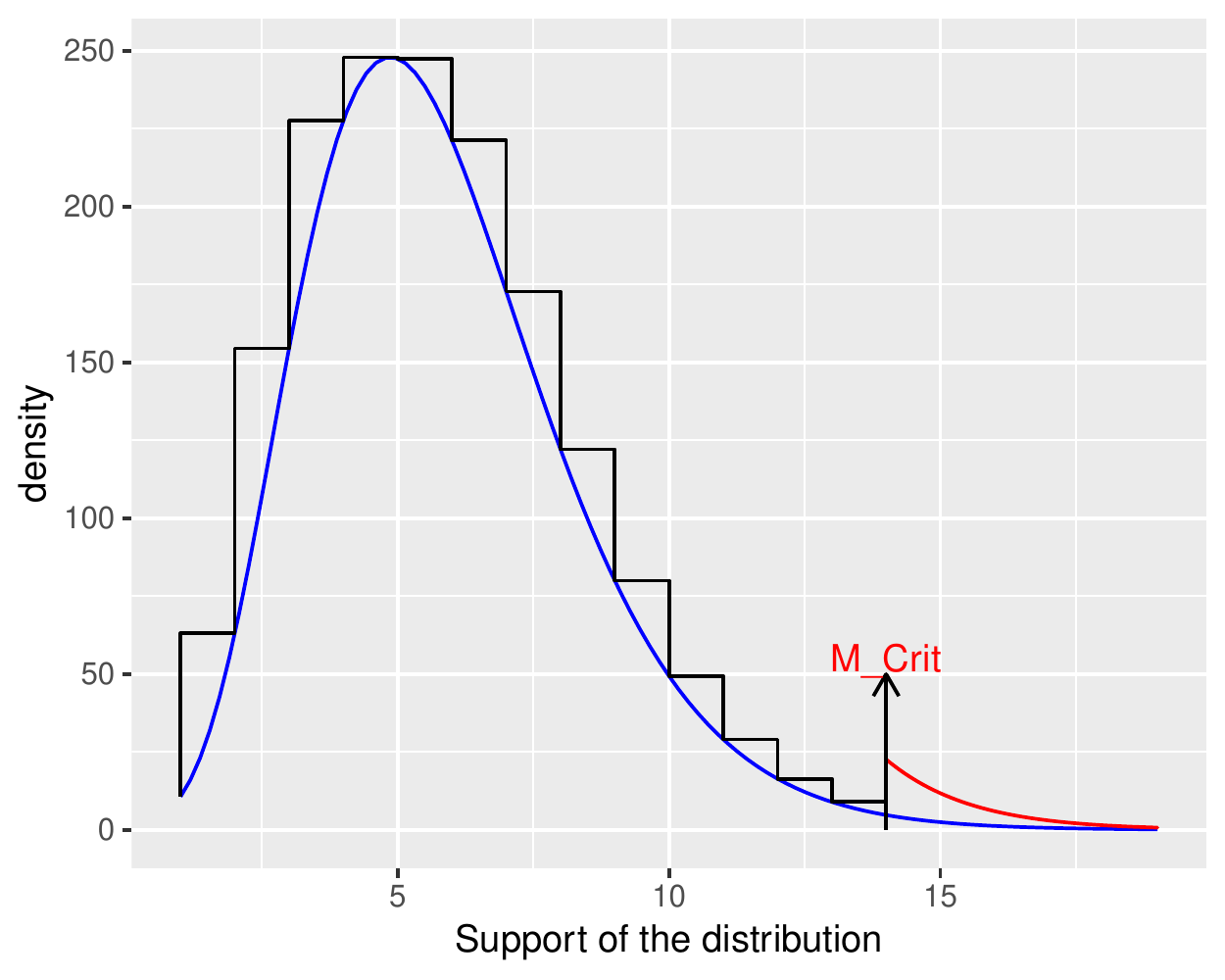} & \includegraphics[scale=0.5]{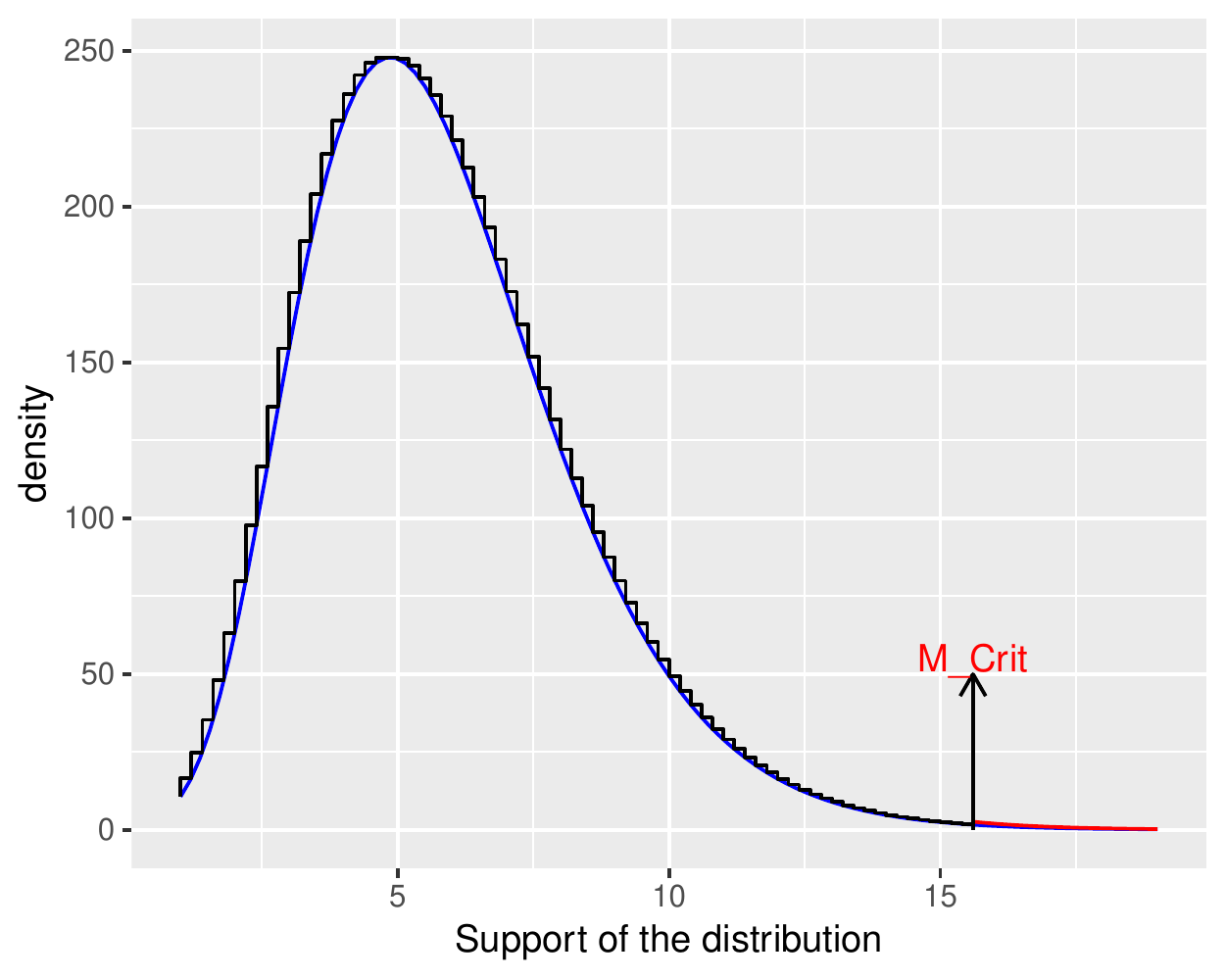}\\
(a)  & (b) 
\end{tabular}
\caption{The blue curves represent $g_{1}$,  the unnormalized density of $CCPD_1^{\star}$ distributions. The black curve and the red curve after $M_{crit}$ constitutes the function $\overline{g}_1$, the proposal density for the accept reject algorithm. The panel(a) displays the key aspects of the densities while panel(b) shows the proposal density when recommended specifications of $M_{crit}$ and $\delta$ are used. }
\label{fig:Diagram}
\end{figure}

We now describe an efficient rejection sampling procedure from the conditional distribution of $(d_1\,\mid\,(d_2, \cdots, d_p))$ when $\bd \sim \IMDY(\cdot ; \nu, \BoEta)$ for some $\nu>0$ and $\max\limits_{1 \leq j \leq p} \eta_j<1$. Here $\BoEta=\left(\eta_1,\ldots, \eta_p \right)$. 
Let $m$ be the mode of the conditional  distribution, $g_1(\cdot) := g(\cdot\,;\,\nu,\BoEta\,\mid\,(d_2, \ldots , d_p))$, of the variable $d_1$ given $(d_2, \ldots , d_p)$ when $\eta_1 > 0$. In case $\eta_1 \leq 0$, we  set $m$ to be $0$.
Using the properties of the conditional distribution described in Theorem~\ref{lem:right_tail_prob_bound}, we compute a critical point $M_{crit}$ such that $P\left(d_{1} > M_{crit} \,\mid \,(d_2, \cdots, d_p), \{{X_j\}}_{j=1}^N \right)< \epsilon$. Here we have chosen $\epsilon=0.0001$. 
{
To construct a proposal density $ \overline{g}_1(x)$, we employ two different strategies, one for the  bounded interval $(0, M_{crit}]$ and the other using Theorem~\ref{lem:right_tail_prob_bound} to tackle the tail, $(M_{crit}, \infty)$, of the support of the conditional posterior distribution of $d_1$. 

The procedure is as follows. Let $\delta = {M_{crit}}/{N_{bin}}$ where $N_{bin}$ is the total number of partitions of the interval $(0, M_{crit}]$. Consider $k=(\left[ {m}/{\delta}\right]+1)$ where $\left[{m}/{\delta} \right]$ denotes the greatest integer less than or equal to $ {m}/{\delta}$. 
Now define the function 
\begin{eqnarray}\label{eq:Defg1Upper}
\overline{g}_1(x) &:=& \sum_{j=1}^{k-1} g_1(j\,\delta)\; \mathbb{I}_{\left( (j-1) \delta, j \delta]\right)} (x)  + g_1(m) \mathbb{I}_{\left( (k-1) \delta, k \delta]\right)} (x) \nonumber\\ 
&&\qquad +\sum_{j=k+1}^{N_{bin}} g_1((j-1)\,\delta)\;  \mathbb{I}_{\left( ( (j-1) \delta, j \delta]  \right)} (x)\nonumber\\ 
&&\qquad + K^{\dagger}_{n,p,M_{crit}} \;d_1^{\nu(n-1)/2} exp(\;-\nu(1-\eta_1)\,d_1)\mathbb{I}_{\left(M_{crit}, \infty)  \right)} (x) ,
\end{eqnarray}
where $K^{\dagger}_{n,p,M_{crit}}$ is as defined in part(d) of Theorem~\ref{lem:right_tail_prob_bound}.

From Theorem~\ref{lem:right_tail_prob_bound} it follows that $ \overline{g}_1(x) \geq g_1(x)$ for all $x>0$ as $g_1(\cdot)$ is a unimodal log-concave function with maxima at $m$. We consider,
$$
q_j = 
\fourpartdef
{{\delta\;g_1(j\delta)}{}}      {1\leq j< \left[\frac{m}{\delta}\right]+1,}
{{\delta\;g_1(m)}{}}      { j= \left[\frac{m}{\delta}\right]+1,}
{{\delta\;g_1((j-1)\delta)}{}} {\left[\frac{m}{\delta}\right]+1< j\leq N_{bin},}
{K^{\dagger}_{n,p,M_{crit}} \frac{\Gamma\left( \frac{(\nu(n-1)+2)}{2},M\nu(1-\eta_1) \right) }{\left[\nu(1-\eta_1)\right]^{\nu(n-1)/2+1}}    } { j= N_{bin}+1,}
$$
where $\Gamma\left( \frac{(\nu(n-1)+2)}{2},M_{crit}\nu(1-\eta_1) \right)$ denotes the upper incomplete gamma function. For the case where  $M_{crit}$ tends to $\infty$ (see  Remark~\ref{remark:KDaggar})  the constant  $K^{\dagger}_{n,p,M_{crit}} $  approaches a finite constant, whereas $\Gamma\left( \frac{(\nu(n-1)+2)}{2},M_{crit}\nu(1-\eta_1) \right)$ monotonically decreases to zero. Therefore,  the positive constant  $q_{_{N_{bin}+1}} $ can be made arbitrary close to zero by choosing a suitably large value for $M_{crit}$ when the value of $ n,p,\nu, \eta_1$ are fixed. 
Note that the quantities $\{q_j\}_{j=1}^{N_{bin}+1}$  may not add up to $1$, therefore  we construct the corresponding set of probabilities,  $\{ p_{j} \}_{j=1}^{N_{bin}+1 } $  where $p_j = {q_j}/{\sum_{j=1}^{N_{bin}+1}q_j}$ for $j = 1, 2, \cdots, N_{bin}+1$.
 The following algorithm lists the steps involved in generating a sample from the distribution corresponding to the kernel $g_1(\cdot)$.
\begin{algorithm}[H]
    \begin{algorithmic}[1]
    \State{ Sample $Z$ from the discrete distribution with the support  $\{ 1, 2, \ldots, (N_{bin}+1)\}$ and corresponding probabilities $\{p_j \}_{j=1}^{N_{bin}+1}$,}
     \If{$Z\leq N_{bin}$ }
        \State { Sample $y \sim \text{Uniform}\left(  (Z-1)\,\delta, Z\delta \right)$,}
        \Else{ Sample $y \sim \text{TruncatedGamma}\left(\text{shape}= \frac{\nu(n-1)+2}{2}, \text{\;rate}= \nu(1-\eta_1), \text {support}=(M_{crit},\infty) \right)$}
        \EndIf
            
        \State { Sample $ U \sim \text{Uniform }(0,1) $,}
        \If{$U\leq \frac{g_1(y)}{\overline{g}_1(y)}$ }
        \State{Accept $y$ as a legitimate sample from $g_{1}(\cdot)$    }
        \Else{ Go to Step 1}
        \EndIf
    \end{algorithmic}
    \caption{Steps for the rejection sampler for $\CCPD^{\star}_j$}
\end{algorithm}

Figure~\ref{fig:Diagram}  shows a typical example of the function $g_1(x)$ and the corresponding $\overline{g}_1(x)$. The blue curve represents the unnormalized density $g_{1}$. The black curve and the red curve after $M_{crit}$ constitutes the function $\overline{g}_1$ ( defined in Equation~\ref{eq:Defg1Upper}). Note that the red curve after the point $M_{crit}$ represents the last term (involving $K^{\dagger}_{n,p,M_{crit}}$) in the summation formula in Equation~\ref{eq:Defg1Upper}. In Figure~\ref{fig:Diagram}(a), the values of $\delta$ and $M_{crit}$ are set such that the key components of $g_1$ and $\overline{g}_1(x)$ are easy to discern. On the other hand, Figure~\ref{fig:Diagram}(b) displays the plot of $\overline{g}_1(x)$ when recommended specification of $M_{crit}$ and $\delta$ are used.
}
%



The  choice of $N_{bin}$ plays a crucial role in the algorithm and is required to be determined before constructing the proposal density for the accept-reject algorithm. Note that $N_{bin}$ and  $\delta$ are  interconnected. If one is specified, the value of the other can be determined. We decide to choose the parameter $\delta$ and compute the corresponding $N_{bin}$. In the case where the concentration parameter is high, a finer partition of the proposal histogram (smaller value of $\delta$) is required to keep the acceptance rate of the algorithm high. Based on our empirical results, we recommend selecting $\delta$ to be of the order of $\frac{1}{\sqrt{\nu}}$.  
The acceptance probability remains stable across different choices of $\nu$ when the value $\delta$ is set accordingly (see Figure~\ref{fig:AcceptenceProb}).
The estimated acceptance probabilities, used in Figure~\ref{fig:AcceptenceProb}, were calculated based on $10000$ Monte Carlo samples for each value of $\nu$ varied from $1$ to $100$. 
The relationship between $N_{bin}$ and $\delta$ and $\nu$ is presented in Table~\ref{Tab:NbinANDdelata}.

  
Finally, successful implementation of the sampling algorithm developed in this subsection requires the computation of $\hyp$, a key step for the computation of $g_1(\cdot)$. In Section~\ref{sec:HYPComputation} we discuss the procedure that we have adopted to compute $\hyp$.
%
%

\begin{table}[h]
\begin{tabular}{|l|l|l|l|}
\hline
$\nu$ & $\delta$ & Estimated Acceptance probability & $N_{bin}$ \\ \hline
1                  & 1                     & 0.95813                          & 42         \\ \hline
1                  & 0.5                   & 0.977517                         & 85         \\ \hline
1                  & 0.333333              & 0.984155                         & 127        \\ \hline
1                  & 0.2                   & 0.988924                         & 212        \\ \hline
1                  & 0.1                   & 0.996314                         & 425        \\ \hline
1                  & 0.05                  & 0.998104                         & 851        \\ \hline
3                  & 0.5                   & 0.952835                         & 27         \\ \hline
3                  & 0.333333              & 0.963206                         & 40         \\ \hline
3                  & 0.2                   & 0.977326                         & 67         \\ \hline
3                  & 0.1                   & 0.988924                         & 135        \\ \hline
3                  & 0.05                  & 0.995124                         & 271        \\ \hline
5                  & 1                     & 0.885818                         & 3          \\ \hline
5                  & 0.5                   & 0.941886                         & 7          \\ \hline
5                  & 0.333333              & 0.960246                         & 10         \\ \hline
5                  & 0.2                   & 0.973994                         & 17         \\ \hline
5                  & 0.1                   & 0.989218                         & 35         \\ \hline
5                  & 0.05                  & 0.993246                         & 71         \\ \hline
\end{tabular}
\caption{Values of the $N_{bin}$, $\delta$ and acceptance probability for algorithm to generate values from  $CCPD_j(\eta, \nu)$ for $\nu=1, 3, 5$.  }
\label{Tab:NbinANDdelata}
\end{table}
%

\subsection{\bf Computation of $\hyp$}\label{sec:HYPComputation}
We first describe an efficient and reliable computational procedure to compute the function  $\hyp$ when the argument matrix $\R$ is of dimension $2 \times 2$. The procedure is relevant to many applications considered in the field  ~\citep{Downs:1971, Downs:1972, Jupp:1979,Jupp:1980, Mardia:1977, Mardia:2007, Mardia:2009, Chikuse:1991:as,Chikuse:1991,Chikuse:1998,Chikuse:2003,Sei:2013, Lin:2017}.  We emphasize that the computational procedure described below is applicable for analyzing data on $\mathcal{V}_{n,2}$ for all $n\geq 2$.

Consider the representation  developed in ~\cite{Muirhead:1975} for the Hypergeometric function of a matrix argument 
\begin{eqnarray}\label{eq:Hyp:Series1}
{}_0F_1\left(\p, \R\right)&=&\sum_{k=0}^{\infty}\frac{\rr_1^{k}\rr_2^{k}}{\left(\p-\frac{1}{2}\right)_k\left(\p\right)_{2k}k!}\; _0F_1\left(\p+2k,{\rr_1+\rr_2}\right),
\end{eqnarray}
where $\R$ is a $2\times 2$ diagonal matrix with diagonal elements $\rr_1>0, \rr_2>0$. From \cite{Butler:2003} (see page 361), it can be seen that, 
\begin{eqnarray}\label{eq:Hyp:Bessel:Connection}
{}_0F_1\left(\p+2k,{\rr}_1+\rr_2\right)=\frac{\Gamma\left(\p+2k\right)}{\left(\sqrt{\rr_1+\rr_2}\right)^{(\p+2k-1)}}I_{\p+2k-1}\left(2\sqrt{\rr_1+\rr_2}\right).
\end{eqnarray}
where $I_{\p+2k-1}(\cdot)$ is the modified Bessel function of the first kind with order $(\p+2k-1)$. Hence from Equation~\ref{eq:Hyp:Series1} and Equation~\ref{eq:Hyp:Bessel:Connection}, we get that 
\begin{eqnarray}\label{eq:Hyp:Series}
{}_0F_1\left(\p, \R\right)
&=& \sum_{k=0}^{\infty}\frac{\rr_1^{k}\rr_2^{k}}{\left(\p-\frac{1}{2}\right)_k\left(\p\right)_{2k}k!}\; \frac{\Gamma\left(\p+2k\right)\;\; I_{\p+2k-1}\left(2\sqrt{\rr_1+\rr_2}\right)}{\left(\sqrt{\rr_1+\rr_2}\right)^{(\p+2k-1)}}\nonumber\\
&=& \sum_{k=0}^{\infty} \TermK_k,
\end{eqnarray}

where $ \TermK_k=\frac{\Gamma(\p-.5)\Gamma(\p)}{\Gamma(\p+k-.5) k!} \frac{(\rr_1 \rr_2)^{k}}{\left(\sqrt{\rr_1+\rr_2}\right)^{(\p+2k-1)}}\;{ I_{\p+2k-1}\left(2\sqrt{\rr_1+\rr_2}\right)}$. Note that 
\begin{eqnarray}
\frac{\TermK_{k+1}}{\TermK_{k}}&=&\frac{\Gamma(\p+k-.5) k!}{\Gamma(\p+k+.5) (k+1)!} \frac{ I_{\p+2k+1}\left(2\sqrt{\rr_1+\rr_2}\right)}{ I_{\p+2k-1}\left(2\sqrt{\rr_1+\rr_2}\right)}\frac{\rr_1 \rr_2}{\left({\rr_1+\rr_2}\right)} \nonumber\\
&\stackrel{}{\leq} & \frac{4\rr_1 \rr_2}   {(2\p+2k-1)(2k+2)(2k+\p)(2k+2\p+1)   },
\end{eqnarray}
 where the last inequality follows from $I_{\nu+1}(x)/I_{\nu}(x)< \frac{x}{2(\nu+1)}$ for $x>0, \nu>-1$ (see page 221 in \cite{Ifantis:1990}).
For fixed values of $\rr_1, \rr_2$ we can find $\M$ such that $A_{\M}\leq \epsilon $ and $M^4 \geq (\rr_1\; \rr_2)/( 4 \epsilon_1)$ for some  $\epsilon_1<\frac{1}{2}$ and  a predetermined error bound $\epsilon$. For such a choice of $\M$, if $k$ is any integer such that $k \geq \M$, then
\begin{eqnarray}
\frac{\TermK_{k+1}}{ \TermK_{k} } &\leq &   \frac{4\rr_1 \rr_2}   {(2c+2k-1)(2k+2)(2k+c)(2k+2c+1)    }\nonumber\\
&\leq &   \frac{4\rr_1 \rr_2}   {(2c+2\M-1)(2\M+2)(2\M+c)(2\M+2c+1)    }\nonumber\\
&\stackrel{}{\leq} &   { \left(\frac{\rr_1 \rr_2}{4M^4} \right)   \left\{ \frac{16 M^4}   {(2c+2\M-1)(2\M+2)(2\M+c)(2\M+2c+1)    }\right\}}\nonumber\\
&\stackrel{}{\leq} &   { \left(\frac{\rr_1 \rr_2}{4M^4} \right)   \left\{ \frac{ M^4}   {(\M+\frac{2c-1}{2})(\M+1)(\M+\frac{c}{2})(\M+\frac{2c+1}{2})    }\right\}}\nonumber\\
&\stackrel{}{\leq}& \epsilon_1,
\end{eqnarray}
where the last inequality follows due to the fact that $  M^4 \leq  (\M+\frac{2c-1}{2})(\M+1)(\M+\frac{c}{2})(\M+\frac{2c+1}{2})    $  as $c> \frac{1}{2}$. Hence from Equation~\ref{eq:Hyp:Series} we get that 
\begin{eqnarray}
\vert {}_0F_1\left(\p, \R\right)- \sum_{k=0}^{\M} \TermK_k \vert = \sum_{k=\M+1}^{\infty} \TermK_{k}\leq \TermK_{\M}\sum_{k=\M+1}^{\infty} \epsilon_1^{k-\M} \leq \frac{\epsilon\;\epsilon_1\;}{1-\epsilon_1}<\epsilon.
\end{eqnarray}
Consequently, for a given value of the matrix $\R$ and an error level $\epsilon$, we can select $\M$  accordingly, so that $ {}_0F_1\left(\p, \R\right)$ is approximated as 
\begin{eqnarray}\label{eq:Hyp:Series:Approximation}
{}_0F_1\left(\p, \R\right)
&\approx& \sum_{k=0}^{\M}\frac{\rr_1^{k}\rr_2^{k}}{\left(\p-\frac{1}{2}\right)_k\left(\p\right)_{2k}k!}\; \frac{\Gamma\left(\p+2k\right)\;\; I_{\p+2k-1}\left(2\sqrt{\rr_1+\rr_2}\right)}{\left(\sqrt{\rr_1+\rr_2}\right)^{(\p+2k-1)}},
\end{eqnarray}

where the error in the approximation is at most $\epsilon$.

In the case when the matrix $\R$ is of dimension $p\times p $ with $p>2$, we rely on the computational technique developed in \citep{Koev:2006}. Development of efficient computational schemes for the hyper geometric function of a matrix argument in general dimension is an active area of research \citep{Gutierrez:2000, Koev:2006, Nagar:2015,Pearson:2017}. In principle, the theoretical framework developed in this article integrated with the general computation scheme specified in \cite{Koev:2006} can handle data on $\mathcal{V}_{n,p}$ for arbitrary integers $n\geq p\geq 2$, but the results from the combined procedure may lack precision as it inherits the limitations of the algorithm in \cite{Koev:2006} ( See  page 835 in \cite{Koev:2006}). In the following remark we  specify the assumptions under which the combined procedure can be applied effectively.  

\begin{rmrk}\label{Remark:Hyp:Koev}
The algorithm developed in \cite{Koev:2006} is a general procedure for computing  $_pF_q (\cdot)$ for arbitrary integers $p,q\geq 0$. Naturally, the algorithm applies to $_0F_1$ which is the object of focus in the current context. Due to its generality, the computational scheme has certain limitations. In particular, it requires appropriate specification of a ``tuning parameter" that can not be determined in an automated manner. 
   However, from an empirical exploration of the procedure, we observed that the corresponding outputs can be quite robust. Particularly, the output was found to stabilize after a certain point (we will call this the ``stabilization point") when the value of the tuning parameter was gradually increased.
   For the case of $p=2$,  if the tuning parameter is specified to be larger than the stabilization point, the output from~\cite{Koev:2006} is very close to the true value, as determined by our arbitrary precision algorithm.
   Extrapolating to $p\geq3$, we presume that the true value of the corresponding hyper geometric function will be close to the output of~\cite{Koev:2006} if the tuning parameter is set larger than the ``stabilization point".
  As the ``stabilization point" is observed to be larger for larger values of $D$,  we can set the value of the tuning parameter to a single pre-specified number for an entire analysis only if we assume that the diagonal elements of the matrix $D$ are bounded above by a prespecified finite number. Under this assumption, we can rely on~\cite{Koev:2006} for the analysis of data on $\mathcal{V}_{n,p}$, $n\geq p\geq 3$. 
In that case, the combination of our theoretical framework and the algorithm for the computation of the hypergeometric function from \cite{Koev:2006} would work effectively for practical applications (see Simulation Section\ref{Sim2}).
   
 
In contrast, the procedure to compute $\hyp$ that we have developed, though targeted towards a specific case, has a theoretical guarantee for a desired level of precision of its output. Since many statistical applications, as mentioned earlier, are about analyzing data on $\mathcal{V}_{n,2}$,
the computation procedure we have designed specifically for $\mathcal{V}_{n,2}$ has its own merit. 
\end{rmrk}

\section{Simulation}
\label{sec:simul_data_app}

To evaluate the performance of the procedure presented in the previous sections, we performed simulation experiments. We considered two different setups. In the first, we analyzed simulated datasets in $\StiefelS$ where we varied $n$ to assess its effect on the posterior estimation efficiency. Here, the value of $p$ was fixed at $2$ and the computation of $\hyp$ developed in Section~\ref{sec:HYPComputation} was utilized. In the second setup, we analyzed data on  $\StiefelS$ to demonstrate the generic applicability of our framework by setting  $p=3$ , $n=5$.  Here, we used the procedure in~\cite{Koev:2006} to calculate the value $\hyp$. 

\subsection{Simulation Setup $(p=2$)}\label{Sim1}
We  present results from experiments with simulated data where we varied the dimension of the Stiefel manifold, $n$, across a range of values. The objective of this simulation study was to see how the error rates varied  with the dimension $n$.  Specifically, we generated $3000$ observations using $\ML$ distribution on $\mathcal{V}_{3,2}$, $\mathcal{V}_{5,2}$, $\mathcal{V}_{10,2}$, and $\mathcal{V}_{15,2}$. These correspond to the Stiefel Manifolds with dimension $[n=3,p=2]$, $[n=5,p=2]$, $[n=10,p=2]$, and $[n=15,p=2]$, respectively. 
We generated $50$ datasets for each simulation setting using the algorithm mentioned in~\cite{Hoff:2009}. In order to generate data for each dataset we fixed the parameters $M$ and $V$ to the canonical orthogonal vectors of appropriate dimension and generated two entries of the parameter $D$ from two independent gamma distributions.

\begin{figure}[!htb]
\centering
\begin{tabular}{cc}
\includegraphics[scale=0.08]{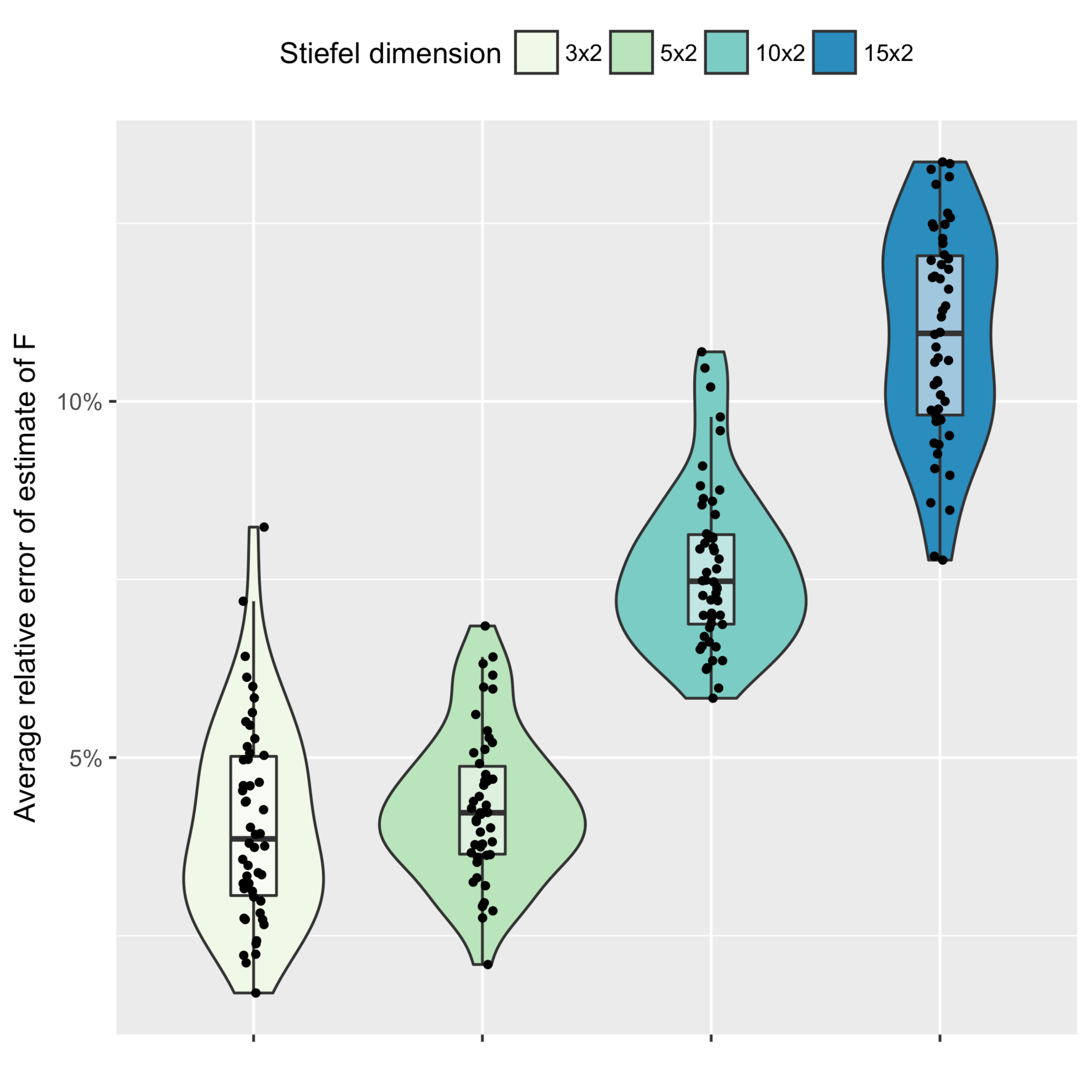} & \includegraphics[scale=0.08]{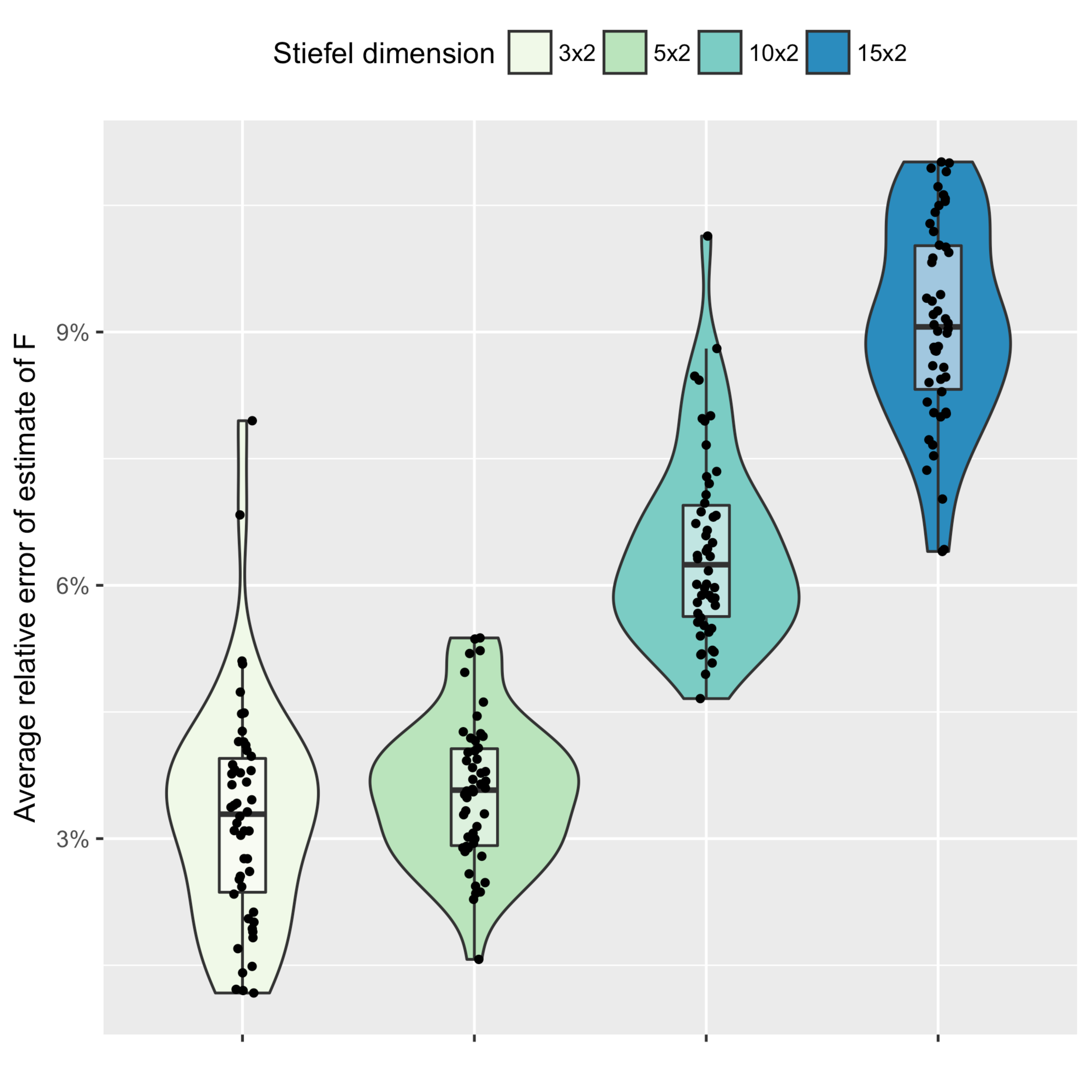}\\
(a) with 2000 data points & (b) with 3000 data points
\end{tabular}
\caption{Relative error of $\widehat{F}$ for matrices with different dimensions }
\label{fig:F_estimate_relative_error}
\end{figure}

We ran posterior inference for each of these datasets using $3000$ MCMC samples with an initial $1000$ samples as burn-in. We used the posterior mean of the parameter $F$ as the point estimate $\widehat{F}$. Finally we assessed our performance by computing the relative error for the estimate of $F_{true} = M_{true}D_{true}V_{true}^T$. We define the relative error as:
\begin{equation*}
\frac{\Vert\widehat{F}-F_{true}\Vert}{\Vert F_{true}\Vert},
\end{equation*}
where $\Vert \cdot \Vert$ denotes the matrix Frobenious norm. Figure~\ref{fig:F_estimate_relative_error} shows the average relative error with the corresponding standard deviation of estimation for $\mathcal{V}_{3,2}$, $\mathcal{V}_{5,2}$, $\mathcal{V}_{10,2}$, and $\mathcal{V}_{15,2}$ for $N=2000$ (panel (a)) and for $N=3000$ (panel (b)). The average relative errors do not seem to exceed $11\%$ and $9\%$ for $N = 2000$ and $3000$, respectively even with the dimension as high as $15$. The error rate tends to increase with higher dimension, i.e., value of $n$. 
Also, we investigated the relationship with the total sample size and found these error rates to decrease with larger sample sizes. For example, the reduction in average relative error rate for $n = 5$ and $N = 2000$ is around $2\%$. Overall, these results demonstrate the robustness of our inference procedure. 

\subsection{Simulation Setup ($p>2$)}\label{Sim2}

\begin{figure}[!htb]
\centering
\includegraphics[scale=0.5]{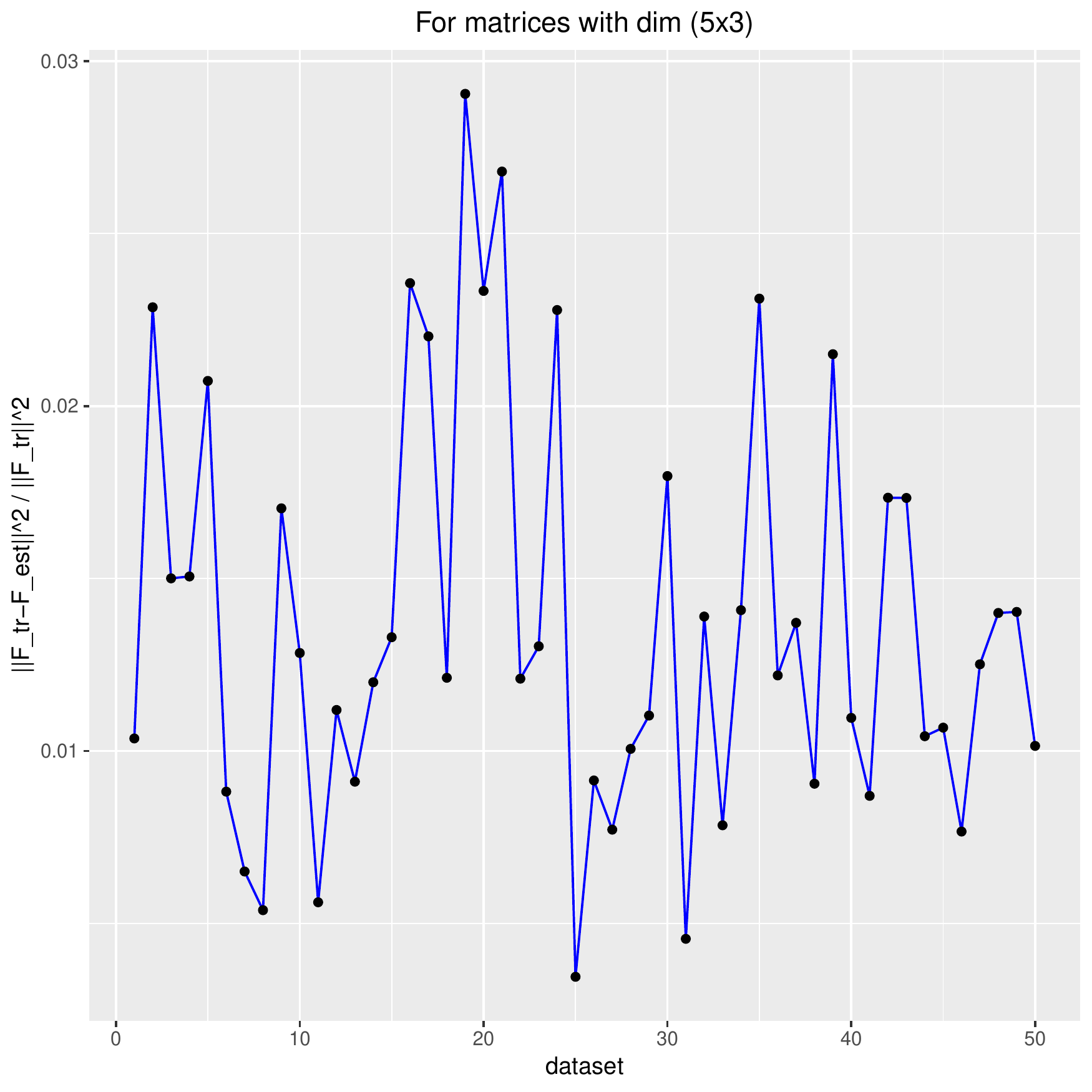} 
\caption{Average relative error for datasets on $\mathcal{V}_{5,3}$}
\label{fig:F_rel_dist_5_3}
\end{figure}

Having demonstrated the efficiency of our method for a range of values of $n$ with $p=2$, we now present an example of a generalized simulation scenario for $p>2$. Here we use the procedure in \cite{Koev:2006} to numerically approximate the value of $\hyp$ where $D$ is a $p\times p$ dimensional matrix with $p>2$ (See Remark~\ref{Remark:Hyp:Koev}). Through the entire simulation we fixed the tuning parameter required in the computation of $\hyp$ to a large prespecified value. Here we give a specific example with $n=5$ and $p=3$. We generated $50$ datasets of $500$ observations each using the $\ML$ distribution with different parameters, on $\mathcal{V}_{5,3}$. We then ran posterior inference for each of these datasets using $1100$ MCMC samples with an initial $100$ sample burn-in. We used the posterior mean of the parameter $F$ as before as the estimate of the true parameter $F$. Using the same metric we computed the average relative error of the estimation (Figure~\ref{fig:F_rel_dist_5_3}). We observed that our sampling algorithm for $d_i$ ($i=1, 2, 3$) runs with a very low rejection rate. As can be seen in Figure~\ref{fig:F_rel_dist_5_3}, the average relative errors do not exceed 3\%, demonstrating the general applicability of our framework beyond $p=2$.

Codes for the algorithms are available at \url{https://github. com/ssra19/Stiefel_Bayes.git}. 

\section{Application}
\label{sec:real_data_app}
Finally, to showcase the methodology developed in this paper, we analyzed the vectorcardiogram dataset discussed in~\cite{Downs:1971}. The dataset contains vectorcardiograms of $56$ boys and $42$ girls aged between $2$ and $19$ years. Individuals in the dataset are partitioned into four groups: groups 1 and  2 consist of boys aged between $2-10$  and  $11-19$ years, while groups 3 and 4 consist of girls aged between $2-10$ and  $11-19$ years. Each sample contains vectorcardiograms acquired using two different measurement systems, the Frank lead system~\citep{Frank:1956, Downs:1971} and the McFee lead system~\citep{Downs:1971}. Here, we restrict ourselves to groups 1 and 3 and measurements acquired using the McFee lead system. For each individual sample, we considered the pair of orthogonal vectors that provides the orientation of the ``QRS loop"~\citep{Downs:1971} in $\mathbb{R}^3$. Each orientation in the sample is defined by a $3 \times 2$ matrix with orthonormal columns, i.e., an element in $\mathcal{V}_{3,2}$. Additional details regarding the measurements, data structures, and data processing can be found in~\cite{Downs:1971}.

\subsection{Parameter estimation}\label{sec:ParameterEstimation}
We modeled the vectorcardiogram dataset using $\ML$ distributions on $\mathcal{V}_{3,2}$. There were $28$ and $17$ observations in groups 1 and 3, respectively. We assumed that each i.i.d observation in group 1 follows a $\ML$ distribution with parameters $M_{group1}, {\bd}_{group1}$ and $V_{group1}$, and likewise, i.i.d observations in group 3 follow a $\ML$ distribution with parameters $M_{group3}, {\bd}_{group3}$ and $V_{group3}$. We used the uniform improper prior for estimation of the parameters related to both groups (see Section~\ref{sec:hyperparameter_singleML}). From Equation~\ref{eq:Improper:Posterior}, we note that the posterior distributions of $(M_{group1}, {\bd}_{group1}$, $V_{group1})$ and $(M_{group3}, {\bd}_{group3},V_{group3})$ given the data are  $$JCPD\left( \cdot \;;\; 28,\; \overline{W}_{group1} \right) \text{ and } JCPD\left( \cdot \; ;\; 17,\; \overline{W}_{group3} \right)  \text{ where }$$
\[\overline{W}_{group1} = \begin{bmatrix}
  0.687 & 0.576 \\
 0.551 &  -0.737 \\
 0.122  &  0.142  
\end{bmatrix}  \text{ and } 
\overline{W}_{group3} = \begin{bmatrix}
 0.682  & 0.585 \\
 0.557 & -0.735 \\
 0.125 &  0.055 
\end{bmatrix}
\]
are the sample means of the observations in groups 1 and 3, respectively. We verified the spectral norm condition in Theorem~\ref{thm:DY_joint_prior} for the posterior distributions to be well defined; we found $\normtwo{\overline{W}_{group1}}=0.946$ and $\normtwo{\overline{W}_{group3}}=0.941$.

Using Theorem~\ref{thm:DY_D_property}, we can infer that the above-mentioned posterior distributions have unique modes. Also from Theorem~\ref{thm:DY_D_property} we can compute the posterior mode and they were

\[\widehat{M}_{group1}=\begin{bmatrix}
-0.650 & 0.733\\
 0.743 &  0.668\\
-0.157  & 0.127 
\end{bmatrix}, 
\widehat{\bd}_{group1}=\begin{bmatrix}
 16.329  \\
  5.953
\end{bmatrix},
\widehat{V}_{group1}=\begin{bmatrix}
 -0.059 & 0.998 \\
 -0.998 & -0.059
\end{bmatrix}.\; 
\]
Similarly, we can compute the posterior mode for the parameters of group 3 (not reported here). To estimate the posterior mean for the parametric functions $$F_{group1}=M_{group1}D_{group1}V_{group1}^T \text{ and } F_{group3}=M_{group3}D_{group3}V_{group3}^T,$$
we ran the MCMC based posterior inference procedure described in Section~\ref{subsubsec:post_comp} to generate MCMC samples from each of the posterior distribution.


For group 1, the posterior mean for the parametric function  $F_{group1}=M_{group1}D_{group1}V_{group1}^T$ was  
$$\doublehat{F}_{group1} = \begin{bmatrix}
5.183 & 9.086    \\
3.583 & -10.996  \\
0.919 & 2.221  
\end{bmatrix}, \text{  } SD(\doublehat{F}_{group1})= \begin{bmatrix}
   
1.527 & 2.354 \\	
1.475 & 2.665 \\	
0.596 & 0.898
\end{bmatrix}, $$ 
where the entries of the matrix $SD(\doublehat{F}_{group1})$ provides the standard deviation for the corresponding entries of $\;\doublehat{F}_{group1}$. From the MCMC samples, we also estimated the posterior density of each entry of $F_{group1}$ and $F_{group3}$. Figure~\ref{fig:FDensityPlot_group1} shows the corresponding density plots.  
The estimates related to group 3 were
$$\doublehat{F}_{group3} = \begin{bmatrix}
  3.249 & 8.547   \\	
3.798 & -10.658 \\	
1.605 & 0.796  	
\end{bmatrix} \text{ and } SD(\doublehat{F}_{group3})= \begin{bmatrix}
   
1.263 & 2.123 \\	
1.359 & 2.624 \\	
0.603 & 0.83 	
\end{bmatrix}. $$ 
	
%

\subsection{Hypothesis testing}
Finally, we conducted a two sample hypothesis test for comparing different data groups on the Stiefel. We have chosen hypothesis testing as one of our demonstrations because a general two sample test that does not rely on asymptotics or on the concentration being very large or very small, has not been reported in the literature for data lying on the Stiefel manifold~\citep{Khatri:1977, Chikuse:2012}. The procedure described here is valid for finite sample sizes and does not require any additional assumptions on the magnitude of the parameters.


We considered the VCG dataset and carried out a test to compare the data group 1 against the data group 3 , i.e.  

 $$ H_{0}: F_{group1}=F_{group3} \;{ \text{vs} }\; H_{A}: F_{group1}\neq F_{group3}.$$

To test the hypotheses in a Bayesian model selection framework, we considered two models ${Model}_0$ and ${Model}_{1}$. In  ${Model}_0$, we assumed $M_{group1}= M_{group3}$, $
\bd_{group1}= \bd_{group3}$,  $V_{group1}=\,V_{group3}$ while in ${Model}_{1}$, we did not impose any structural dependencies between the parameters. We assumed the prior odds between the models to be $1$ and computed the Bayes factor 
 $$ B_{0,1}=\frac{P(Data \mid {Model}_0)}{P(Data \mid {Model}_1)},$$
where  $Data$ denotes the combined data from both groups.  Since an analytic form for the Bayes factor is not available in this case, we used an MCMC based sampling technique to estimate the Bayes factor. We used the empirical prior (see Section~\ref{sec:hyperparameter_singleML}) with the choice of prior concentration set at $1$ percentage of the corresponding sample size.
 We followed the procedure described in Section~\ref{subsubsec:post_comp} to generate MCMC  samples from each of the required posterior distribution.
 We used the harmonic mean estimator (HME) \citep{Newton:Raftery:1994} to estimate the marginal likelihoods required for computing the Bayes factor. It is well known that the HME may not perform well when using improper priors. Consequently, unlike in Section~\ref{sec:ParameterEstimation} where we focus on the parameter estimation, we use an informative prior for this part of the analysis. We observed that the HME estimator is stable for the current context. The estimate of $log(B_{01})$ was $51.994$. Hence, we conclude that there is not enough evidence to favor ${Model}_{1}$ over ${Model}_{0}$.


%
%
%
%
%
%
%

\section{Discussion and Future Directions}
\label{sec:disc}
In this article, we have formulated a comprehensive Bayesian framework for analyzing data drawn from a $\ML$ distribution. We constructed two flexible classes of distributions, $\CCPD$ and $\JCPD$, which can be used for constructing conjugate priors for the $\ML$ distribution. We investigated the priors in considerable detail to build insights into their nature, and  to identify interpretations for their hyper-parameter settings. Finally, we explored the features of the resulting posterior distributions and developed efficient computational procedures for posterior inference. An immediate extension would be to expand the framework to mixtures of $\ML$ distributions, with applications to  clustering of data on the Stiefel manifold.

On a related note, we observed that the tractability of the set of procedures proposed in this article depends crucially on one's capacity to compute the hypergeometric function $ \hypinlineF$ as a function the matrix $F$. We were naturally led to a modified representation of $\hypinline$(see Section~\ref{sec:stiefel_distr}) as a function of a vector argument $\bd$. We explored several  properties of the function $\hypinline$, that are applicable to research areas far beyond the particular problem of interest in this article. As a special note, we should highlight that we designed a  tractable procedure to compute the hypergeometric function of a $n\times 2$ dimensional matrix argument. There are many applications in the literature~\citep{Mardia:1977,Jupp:1979, Chikuse:1998,Chikuse:2003, Lin:2017} where the mentioned computational procedure of $\hyp$ can make a significant impact. As such, the manner in which we have approached this computation is entirely novel in this area of research and the procedure is scalable to ``high-dimensional" data, such as in diffusion tensor imaging. In the near future, we plan to further explore useful analytical properties of the hypergeometric function, and extend our procedure to build reliable computational techniques for the hyper-geometric function where the dimension of the matrix argument is $n\times p$ with $p\geq 3$.  

Finally, there is scope for extending the newly proposed family of prior distributions to a larger class of  Bayesian models involving more general densities on manifolds. The properties of the prior and posterior discovered can also be seamlessly generalized. The coming together of state-of-the-art Bayesian methods incorporating topological properties of the underlying space is a rich  area or research interest. We plan to continue to explore this direction and contribute to a Bayesian methodological development on general analytic manifolds.





\bibliographystyle{plainnat}
\bibliography{ref_param_dti}

\end{document}